\documentclass[a4paper,11pt,reqno,twosided]{amsart} 
\usepackage{amsmath,amssymb,amsthm,latexsym}
\usepackage[numbers,sort&compress,longnamesfirst]{natbib}

\newcommand{\embed}{\hookrightarrow}  

\newcommand{\charak}{\chi}
\newcommand{\dist}{\text{dist}}

\newcommand{\xpfeil}{\xrightarrow}  
  
\newcommand{\rand}{\partial} 
\newcommand{\where}{\,|\,}

\newcommand{\schwach}{\rightharpoonup}

\newcommand{\supp}{\mbox{supp}}  
\newcommand{\aequi}{\Longleftrightarrow}  
  
\newcommand{\laplace}{\Delta}

\newcommand{\nz}{{\mathbb N}}

\newcommand{\rz}{{\mathbb R}}

\newcommand{\eps}{\varepsilon}  
\renewcommand{\phi}{\varphi}

\theoremstyle{plain}
\newtheorem{theorem}{Theorem}[section]  
\newtheorem{corollary}[theorem]{Corollary}  
\newtheorem{lemma}[theorem]{Lemma}  
\newtheorem{proposition}[theorem]{Proposition}
\newtheorem*{theorem*}{Theorem}
\newtheorem{remark}[theorem]{Remark}

\theoremstyle{remark}

\newtheoremstyle{citing}
  {3pt}
  {3pt}
  {\itshape}
  {}
  {\bfseries}
  {.}
  {.5em}
  {\thmnote{#3}}

\theoremstyle{citing}

\numberwithin{equation}{section}

\begin{document}
 
\title[Indefinite elliptic problems]
{Existence and Nonexistence of Positive Solutions \\
of Indefinite Elliptic Problems in $\rz^N$}

\author{Matthias Schneider}

\thanks{Research supported by a S.I.S.S.A. postdoctoral fellowship.}
\address{Scuola Internazionale di Studi Avanzati\\
S.I.S.S.A.\\
Via Beirut 2-4\\
34014 Trieste, Italy}
\email{schneid@sissa.it} 

\date{May 2002}  
\begin{abstract}
Our purpose is to find  positive solutions $u \in D^{1,2}(\rz^N)$ of the semilinear elliptic problem
$-\laplace u - \lambda V(x) u = h(x) u^{p-1}$ for $2<p$. The functions $V$ and $h$ may have an
indefinite sign and the linearized operator need not to have a first (principal) eigenvalue,
e.g. we allow $V\equiv 1$.
We give precise existence and nonexistence criteria, which depend on $\lambda$
and on the growth of $h^{-}$ and $h^{+}/V^+$. Existence theorems  
are obtained by constrained minimization. The mountain pass theorem leads to a second solution.
\end{abstract}

\maketitle

\section{Introduction}
We are interested in finding positive solutions of
\addtocounter{equation}{1}
\begin{align}
\label{eq:1}
\tag*{$(\theequation)_\lambda$}
\begin{split}
-\laplace u(x) -\lambda V(x) u(x) + h(x) u(x)^{p-1}= 0 \text{ in } \rz^N, \\
0 \lneqq u \in D^{1,2}(\rz^N)\cap L^2(\rz^N, |V|) \cap L^p(\rz^N,|h|),   
\end{split}
\end{align}     
where $\lambda > 0$ is a real nonnegative parameter. The functions $h$ and $V$ may change sign.
We denote by ${D^{1,2}(\rz^N)}$ the closure of ${C^{\infty}_c(\rz^N)}$ with respect to the norm
${(\int |\nabla u|^{2})^{\frac{1}{2}}}$ in $L^{2^*}$, where we set $2^*:= 2N/(N-2)$. Moreover, for a
given measurable nonnegative function $k$ and $\Omega \subset \rz^N$ we will denote by
$L^q(\Omega,k)$ the space of measurable functions $u$ satisfying
\[ \|u\|_{L^p(\Omega,k)}^p := \int_\Omega k(x) |u|^p  <\infty,\]
where functions, which are equal $k(x) dx$-almost everywhere, are
identified as usual.
For any function $f: \rz^N \to \rz$ we abbreviate its positive part, $\max(0,f)$, by $f^+$
and its negative part, $(f^+-f)$, by $f^-$. Our basic requirements are
\begin{align} \label{eq:2}
N \ge 3,\; p>2 \text{ and } h,\,V \in C(\rz^N):\; V^+\not\equiv 0 \not\equiv h^+,\\  
\label{eq:3}
D^{1,2}(\rz^N) \text{ is compactly embedded in } L^p(\rz^N,h^-),\\
\label{eq:4}
D^{1,2}(\rz^N) \cap L^p(\rz^N, h^+) \text{ is compactly embedded in } L^2(\rz^N, V^+).
\end{align}
We first state some sufficient conditions for the validity of (\ref{eq:3}) and (\ref{eq:4}). Obviously
(\ref{eq:3}) holds whenever $h^-\equiv 0$. 
Moreover, (\ref{eq:3}) is valid under each of the following two assumptions 
(see \cite[Cor. 2.2]{Schn01}, for a necessary and sufficient condition see \cite[Th. 2.1]{Schn01})
\begin{align}
\label{eq:46}
p<2^* \text{  and  }  h^- \in L^{\frac{2^*}{2^*-p}}(\rz^N),\\
\label{eq:47} 
p<2^*,\, h^- \in L^\infty(\rz^N) \text{ and } \limsup_{|x|\to \infty} h^-(x)
  |x|^{\frac{N-2}{2}p-N}=0.
\end{align}
Condition (\ref{eq:4}) holds under each of the following two
assumptions (see \cite[Th. 2.3]{Schn01}),
\begin{align}
\label{knuutz} 
\exists R>0:\: h(x)>0 \; \forall |x|\ge R\text{ and } 
\int_{\rz^N\backslash B_R(0)} V^+ \left(\frac{V^+}{h}\right)^{\frac{2}{p-2}}< \infty,\\
\label{eq:48} 
V^+\in L^{N/2}(\rz^N). 
\end{align}
Note that (\ref{eq:48}) is sufficient for $D^{1,2}(\rz^N)$ to be
compactly embedded in $L^2(\rz^N,V^+)$, hence (\ref{eq:4}) is
satisfied for any $h^+$ under the assumption (\ref{eq:48}).  
\subsection*{Existence results.}
Equation \ref{eq:1} with $V\equiv 1$ was investigated on bounded domains $\Omega$ with various
boundary conditions using bifurcation theory and the method of sub- super-solutions. We mention the
work of \citet{AmaLogo98,Logo00,Logo97,Ouy91}. They used an abstract result of \citet{CraRab71} and
showed that $(\lambda_1(\Omega),0)$ is a bifurcation point for positive solutions of \ref{eq:1},
where $\lambda_1(\Omega)$ is the first eigenvalue of the Laplacian in $\Omega$ with corresponding
boundary conditions, i.e. there is an open neighborhood $U$ of $(\lambda_1,0)$ such that every
nontrivial solution pair $(\lambda,u)$ in $U$ is described by a continuous curve $(\lambda(s),u(s))$
emanating from $(\lambda_1(\Omega),0)$. It was shown that the solution curve is increasing in
$\lambda$ near $\lambda=\lambda_1(\Omega)$ (supercritical) if and only if
\begin{align}
\label{eq:5}
\int_\Omega h(x) e_1(\Omega)^p >0,  
\end{align}
where $e_1(\Omega)$ denotes the positive eigenfunction belonging to the first eigenvalue
$\lambda_1(\Omega)$. \citet{Ouy91} showed for compact Riemannian manifolds that (\ref{eq:5}) is also
necessary for the existence of positive solutions to \ref{eq:1} for $\lambda>\lambda_1(\Omega)$. The
interesting point in this situation is that the problem becomes affected by the positive part of
$h$, $h^+:= \max(0,h)$, which pushes up the spectrum of the
linearized problem and possibly creates a ``ground state'' solution.\\
\citet{AlaTar96,AlaTar93} and \citet{BerCapNir95,BerCapNir94} studied \ref{eq:1} with $V\equiv 1$ on
bounded domains $\Omega$. They used variational methods like constrained minimization and a
variational formulation of Perron's sub- super-solution method. Apart form various existence
results, (\ref{eq:5}) was derived as a necessary and sufficient conditions for the existence of
positive solutions to \ref{eq:1} for $\lambda>\lambda_1(\Omega)$
with Neumann and Dirichlet boundary conditions.\\
To deal with unbounded domains $\Omega$ and general functions $V$ we have to adjust the value
$\lambda_1(\Omega)$. To this end we define for measurable $\Omega \subset \rz^N$ and $V$
\begin{align}
\notag
D^{1,2}(\Omega)&:= \{u \in D^{1,2}(\rz^N) \where u=0 \text{ a.e. on
  }\rz^N\backslash \Omega\}, \\ 
\label{eq:6}
\lambda_1(\Omega,V)&:= \inf\{\|\nabla u\|_2^{2}\where u \in D^{1,2}(\Omega) \cap L^2(\rz^N,|V|),
\;\int V(x) |u|^2 = 1\}.
\end{align}
Note that for bounded domains $\Omega$ and $V \equiv 1$ the value $\lambda_1(\Omega,V)$ coincides
with $\lambda_1(\Omega)$, the usual first eigenvalue of the Dirichlet Laplacian in
$\Omega$. If $\lambda_1(\Omega,V)$ is attained, we will call $\lambda_1(\Omega,V)$ a
principal eigenvalue and a corresponding nonnegative minimizer a principal eigenfunction, denoted by
$e_1(\Omega,V)$, which satisfies due to its variational characterization 
\[\int \nabla e_1(\Omega,V)\nabla \phi = \lambda_1(\Omega,V) \int V(x) e_1(\Omega,V) \phi \quad \forall
\phi \in D^{1,2}(\Omega)\cap L^2(\rz^N,|V|).\]  
Harnack's inequality \cite[C.1]{Sim82} shows that if $\Omega$ is a domain and $V$ is regular enough,
e.g. (\ref{eq:2}), then the principal eigenvalue $\lambda_1(\Omega,V)$, if it exists, is simple
and the corresponding eigenfunction $e_1(\Omega,V)$ is positive inside $\Omega$. 
If $V^+ \in L^{N/2}(\Omega)$ then $\lambda_1(\Omega,V)$ is
attained, since $D^{1,2}(\Omega)$ is compactly embedded in $L^2(\Omega,V^+)$ in this case (see
\cite[Cor. 2.4]{Schn01}).\\ 
Positive solutions of \ref{eq:1} have up to now only been obtained under the assumption  
that the weight $V^+$ decreases fast enough at infinity to ensure that ${\lambda_1(\rz^N,V)}$
is a principal eigenvalue with positive 
eigenfunction $e_1(\rz^N,V)$ (see \cite{DraHua99,BroSta98,Sta98}). 
\citet{CinGam96} proved under the assumptions $V\in L^{N/2}(\rz^N)\cap
L^\infty(\rz^N)$, $2<p<2^*$ and $h \in
L^{2^*/(2^*-p)}(\rz^N)$
\begin{align}
\text{\ref{eq:1} is solvable in a right neighborhood of }
\lambda_1(\rz^N,V) \text{ if } \int h(x) e_1(\rz^N,V)^{p} >0.
\label{eq:7}   
\end{align}
\citet{CosTeh01} showed
that (\ref{eq:7}) remains valid under the assumptions  
\begin{align*}
2<p<2^*,\,V\in L^{N/2}(\rz^N)\cap L^\alpha(\rz^N) \text{ for some }\alpha>N/2 \text{ and }
\lim_{|x|\to \infty}h(x)=h_\infty>0.   
\end{align*}
Our purpose is to give existence and nonexistence criteria for positive solutions
of \ref{eq:1} in cases where the linearized operator need not to have a first (principal) eigenvalue, 
e.g. for $V \equiv 1$ and $\Omega = \rz^N$. 
We emphasize that our approach yields new results
even when $V^+ \in L^{N/2}(\rz^N)$ and a principal eigenfunction $e_1(\rz^N,V)$ exists
since we do not impose any growth condition on $h^+$ and $V^-$ in this case.\\
We shall find solutions of \ref{eq:1} as critical points
of $I_\lambda \in C^2(E,\rz$), defined by
\[I_\lambda(u) := \frac{1}{2} \|\nabla u\|_2^{2} - \frac{\lambda}{2}
\int V(x) |u|^2 + \frac{1}{p} \int h(x) |u|^{p},\]
where we denote by $E$ the Banach space $D^{1,2}(\rz^N)\cap L^2(\rz^N,|V|) \cap L^p(\rz^N,|h|)$
equipped with the norm $\|u\|_E := \|\nabla u\|_2+
\|u\|_{L^2(\rz^N,|V|)}+\|u\|_{L^p(\rz^N,|h|)}$. Under the assumptions (\ref{eq:3}) and (\ref{eq:4})
the space $E$ coincides with $D^{1,2}(\rz^N)\cap L^2(\rz^N,V^-) \cap L^p(\rz^N,h^+)$ and we may
replace its norm by the equivalent norm  $\|\nabla u\|_2+
\|u\|_{L^2(\rz^N,V^-)}+\|u\|_{L^p(\rz^N,h^+)}$.\\  
For $\lambda>\lambda_1(\rz^N,V)$ and $h^{-} \not\equiv 0$ the energy functional $I_\lambda$ is
neither bounded above, nor below and zero is not a local 
minimum of the energy. To find nontrivial solutions we use a constrained minimization method and
apply the mountain pass theorem to a local minimizer of the energy.
Conditions (\ref{eq:3}) and (\ref{eq:4}) prevent the possible lack of compactness at infinity.
Note that in our framework $I''_\lambda(u)$ need not to be a Fredholm operator from $E$ to
$E'$, see Remark \ref{sec-existence-rem1} for
a precise discussion, which leads to additional technical difficulties.\\ 
Motivated by the work of \citet{Ouy91} we define for measurable ${\Omega \subset \rz^{N}}$
\begin{align*}
\lambda_1(\Omega,V,h) := \inf\{\|\nabla u\|_2^{2}\,\where\, &u \in D^{1,2}(\Omega)\cap
 L^2(\rz^N,|V|)\cap L^p(\rz^N,|h|),\\
&\int V(x) |u|^2=1, \, \int h(x) |u|^{p}\le 0\}.  
\end{align*}
We set $\lambda_1(\Omega,V,h)=+\infty$, if the infimum is taken over an empty set.
Assumptions (\ref{eq:2})-(\ref{eq:4}) imply that $\lambda_1(\rz^N,V,h)$ is
attained or equals $+\infty$ (see Lemma \ref{sec:behav} below).\\  
We can now formulate our main existence result. 
\begin{theorem}
\label{sec:intro:t1}
Suppose (\ref{eq:2})-(\ref{eq:4}) hold. Then \ref{eq:1} has a positive solution $u \in C^1(\rz^N)$ in
$E$ for all $\lambda$ such that $\lambda_1(\rz^N,V)<\lambda<\lambda_1(\rz^N,V,h)$. The solution
$u$ is a local minimum of $I_\lambda$. If, moreover, $h$ is nonnegative then $u$ is the unique
positive solution in $E$.  
\end{theorem}
If (\ref{eq:2}) holds, $V^+ \in L^{N/2}(\rz^N)$ and $h^- \in
L^{2^*/(2^*-p)}(\rz^N)$ then $e_1(\rz^N,V)$ exists, $\lambda_1(\rz^N,V)$ is
simple and $\lambda_1(\rz^N,V,h)$ is attained or equals $+\infty$. Consequently 
\[\lambda_1(\rz^N,V,h)>\lambda_1(\rz^N,V) \text{ if and only if }
e_1(\rz^N,V) \not \in \big\{ u \in L^p(\rz^N,|h|) \where \int h(x)|u|^p \le 0\big\}. \]  
Hence Theorem \ref{sec:intro:t1} includes the existence results in
\cite{CinGam96,CosTeh01} mentioned above. 
From our proof it is easy to see that Theorem \ref{sec:intro:t1} continues to hold for arbitrary
domains $\Omega \subset \rz^N$ and that for a bounded domain $\Omega$, $V\equiv 1$ and $h \in
L^{\infty}(\Omega)$ the condition $\lambda_1(\Omega,1,h)>\lambda_1(\Omega,1)$ is
equivalent to (\ref{eq:5}).\\ 
To give our next existence result we need to introduce the set $\Omega^{-0}$, defined by
\[\Omega^{-0}:= \{x \in \rz^N \where h(x) \le 0\}.\]
If $h$ is nonnegative then we have by definition $\lambda_1(\rz^N,V,h)=\lambda_1(\Omega^{-0},V)$ and we show 
\begin{theorem}
\label{sec:intro:t2}
Suppose (\ref{eq:2})-(\ref{eq:4}) and $h^- \equiv 0$. Then
$\lambda_1(\rz^N,V)<\lambda_1(\Omega^{-0},V)$ and \ref{eq:1} has a unique
positive solution $u_{\lambda}$ in $E$ for
every $\lambda \in J:= \left]\lambda_1(\rz^N,V),\lambda_1(\Omega^{-0},V)\right[$.\\ 
The map $\lambda \mapsto u_\lambda$ belongs to $C^0(J,E)$ and satisfies 
\begin{align*}
\lim_{\lambda \to \lambda_1(\rz^N,V)}u_\lambda =0 \text{, } \lim_{\lambda \to
  \lambda_1(\Omega^{-0},V)}\|u_\lambda\|_E = \infty.
\end{align*}
Furthermore, \ref{eq:1} does not admit any positive solution in $E$ if
$\lambda \ge \lambda_1(\Omega^{-0},V)$. 
\end{theorem}
The continuity of $\lambda \mapsto u_\lambda$ is proven merely by hands since we cannot apply
the implicit function 
theorem, which is mainly due to the fact that $I''_\lambda(u_\lambda)$ fails to be a Fredholm
operator (see Lemma \ref{sec-existence-fred} and Remark \ref{sec-existence-rem1}).\\  
If $V$ is nonnegative then the blow-up behavior of $u_{\lambda}$ as $\lambda \to
\lambda_1(\Omega^{-0},V)$ is given by 
\begin{theorem}
\label{sec:intro:t2:wachs}
Under the assumptions of Theorem \ref{sec:intro:t2}, if moreover $V$ is nonnegative,
$\Omega^{-0}$ is bounded and $\lambda_1(\Omega^{-0},V)<\infty$, then we
obtain the pointwise growth estimate 
\begin{align*}
u_\lambda(x)\ge \frac{(\lambda-\mu)C}{\lambda_1(\Omega^{-0},V)-\lambda} \;\; e_1(\Omega^{-0},V)(x)+ u_\mu(x),
\end{align*}
for any $\mu \in J$, any principal eigenfunction
$e_1(\Omega^{-0},V)$ and some $C=C(e_1,\mu)>0$.   
\end{theorem}
The analogous result for bounded domains was obtained in \cite{GarLogoetal98} using that
$u_\lambda$ is differentiable as a function of $\lambda$. We are able to extend their result
replacing $(d/d\lambda) u_\lambda$ by corresponding difference quotients.\\ 
If $h^- \not\equiv 0$, it is possible to find a second positive solution thanks to the mountain
pass theorem \cite{AmbRab73}. We prove 
\begin{theorem}
\label{sec:intro:t3}
Suppose (\ref{eq:2})-(\ref{eq:4}), $h^- \not\equiv 0$ and $\lambda_1(\rz^N,V)< \lambda_1(\rz^N,V,h)$.
Then for every $\lambda \in J^-:=]\lambda_1(\rz^N,V),\lambda_1(\rz^N,V,h)]$ equation
\ref{eq:1} admits an ordered pair of positive solutions $u_\lambda<v_\lambda$ such that 
\begin{itemize}
\item[] $I_{\lambda}(u_\lambda)<0$ and $(I_\lambda''(u_\lambda)\phi,\phi)\ge 0$ for all 
$\phi \in E$,
\item[] $I_\lambda(v_\lambda)>0$ and $(I_\lambda''(v_\lambda)v_\lambda,v_\lambda)<0$ for all $\lambda
  \in J^-\backslash\{\lambda_1(\rz^N,V,h)\}$,
\item[] $I_{\lambda}(v_\lambda) =0 =(I_\lambda''(v_\lambda)v_\lambda,v_\lambda)$ if
  $\lambda=\lambda_1(\rz^N,V,h)$.  
\end{itemize}
\end{theorem}
Consider the quantity $\lambda^*$, defined by
\begin{align*}
\lambda^* := \sup \left\{ \lambda \where \lambda> \lambda_1(\rz^N,V) 
\text{ and \ref{eq:1} admits a solution in } E \right\}.
\end{align*}
Theorem \ref{sec:intro:t2} shows that $\lambda^*=\lambda_1(\Omega^{-0},V)=\lambda_1(\rz^N,V,h)$
if $h$ is nonnegative. \citet{Ouy91} discussed \ref{eq:1} on compact
Riemannian manifolds $M$ with $V\equiv 1$ and showed that if $h^- \not\equiv 0$ and
$\lambda_1(M)<\lambda_1(M,1,h)$ then 
$\lambda^*>\lambda_1(M,1,h)$ and \ref{eq:1} admits a unique positive solution in $E$ for
$\lambda=\lambda^*$ and at most two ordered positive solution for
$\lambda_1(M)<\lambda<\lambda^*$. The proof consists of a delicate bifurcation analysis which does
not work under our assumptions. A result in the spirit of \cite{Ouy91} is given by the following theorem.
\begin{theorem}
\label{sec:intro:t4}
Suppose (\ref{eq:2})-(\ref{eq:4}), $V \ge 0$, $h^- \not\equiv 0$, $\lambda_1(\rz^N,V)<\lambda_1(\rz^N,V,h)$
and 
\[0<\int h^-(x) e_1(\Omega^{-0},V)^p\]
for some principal eigenfunction $e_1(\Omega^{-0},V)$. Then
\begin{itemize}
\item[(i)] If $\lambda \in ]\lambda_1(\rz^N,V),\lambda^*]$, then \ref{eq:1} has a positive solution
$u_\lambda \in E$ such that 
\[I_{\lambda}(u_\lambda)<0 \text{ and } (I''(u_\lambda)\varphi,\varphi)\ge 0 \text{ for all } 
\phi \in E,\] 
\item[(ii)] $\lambda^*<\lambda_1(\Omega^{-0},V)$ and \ref{eq:1} does not admit any positive solution in $E$
if $\lambda>\lambda^*$.
\end{itemize}
\end{theorem}
Obviously we have $\lambda_1(\rz^N,V,h)\le \lambda^*$. To see that $\lambda^*$ could be arbitrary close to
$\lambda_1(\Omega^{-0},V)$ we
discuss the dependence of $\lambda_1(\rz^N,V,h)$ on $h^+$ showing
\begin{theorem}
\label{sec:intro:l1}
Suppose (\ref{eq:2})-(\ref{eq:4}) hold and $\lambda_1(\rz^{N},V,h)$ is finite. We define for $\mu \in
\rz$ the function $h_\mu := \mu h^{+}-h^{-}$. Then 
\begin{align}
\lambda_1(\rz^N,V,h_\mu) \to \lambda_1(\Omega^{-0},V) \text{  as  }  \mu \to \infty.
\end{align}
Moreover, if $\lambda_1(\Omega^{-0},V)< \infty$, then $\lambda_1(\Omega^{-0},V)$ is attained.
\end{theorem}

Theorem \ref{sec:intro:t1} extends to rather general nonlinearities
(see \cite{SchnDis01}), because the solution is a local
minimizer of the variational functional and may be found whenever the
nonlinearity is weakly lower semicontinuous.
Whereas the proof of the Palais-Smale condition becomes a delicate
issue in presence of a general sign changing nonlinearity (see
\cite{RamTerTro98,AlaTar93}). A version of Theorem \ref{sec:intro:t4}
for special classes of nonlinearities, e.g. $h(x)g(u)$, is given
in \cite{SchnDis01}.

\subsection*{Nonexistence results.}
Roughly speaking there are two obstructions to the existence of (weak) solutions of \ref{eq:1} for
$\lambda>\lambda_1(\rz^N,V)$. The
quotient $h^+/V^+$ has to grow at infinity to ensure the compactness of the inclusion of
$D^{1,2}(\rz^N) \cap L^p(\rz^N,h^+)$ in $L^2(\rz^N,V^+)$ and the parameter $\lambda$ must not exceed
$\lambda_1(\Omega^{-0},V)$. We will prove
\begin{theorem}
\label{sec:intro:t5}
Suppose (\ref{eq:2}). Then \ref{eq:1} does not admit any positive weak solution in $E$
for $\lambda > \lambda_1(\Omega^{-0},V)$.\\
Furthermore, if $\lambda_1(\Omega^{-0},V)<\infty$ 
is attained and $h \not\equiv 0$, 
then there is no positive weak solution of \ref{eq:1} in $E$ for
$\lambda=\lambda_1(\Omega^{-0},V)$.    
\end{theorem}  
If $h$ is nonnegative then $\lambda_1(\rz^N,V,h)= \lambda_1(\Omega^{-0},V)$ and in this case our
results for $\lambda$ are optimal. Moreover, the assumption $h^+ \not \equiv 0$ is necessary because
otherwise $\Omega^{-0}= \rz^N$ and no solution exist for $\lambda > \lambda_1(\rz^N,V)$.\\
In some cases the existence of a positive solution to \ref{eq:1} in $E$ for
$\lambda>\lambda_1(\rz^N,V)$ implies the compactness of the embbeding in (\ref{eq:4}).
To our knowledge such a connection has never been observed in the literature
before. Besides such rather abstract results
we show that the existence of positive solutions leads to
integrability conditions on $V$ and $h$ regaining (\ref{knuutz}) as a 
necessary condition. A related result in dimension $N=1$ can
be found in \cite{kuepp79}. We state our main result, which will be proven below.
\begin{theorem}
\label{sec:intro:t8}
Suppose (\ref{eq:2})-(\ref{eq:3}) and
\begin{align}
\label{eq:28}
0 < \liminf_{|x| \to \infty} \frac{V(x)}{|x|^\alpha} \le \limsup_{|x| \to \infty}
\frac{V(x)}{|x|^\alpha}<\infty \text{ for some }\alpha >-2,\\
\label{eq:40}
0 < \liminf_{|x| \to \infty} \frac{h(x)}{|x|^\beta} \le \limsup_{|x| \to \infty}
\frac{h(x)}{|x|^\beta}<\infty \text{ for some } \beta \in \rz.  
\end{align}
Then $\lambda_1(\rz^N,V)=0$ and (\ref{eq:4}) is equivalent to each of the following two
conditions 
\begin{align}
\label{eq:31}
\lambda_1(\rz^N,V,h)>0 \text{ and \ref{eq:1} is solvable for
  all $\lambda$ such that }0<\lambda<\lambda_1(\rz^N,V,h),\\
\label{eq:12}
\text{\ref{eq:1} is solvable for some } \lambda>0.
\end{align}
If, moreover, $(N-4)\le \alpha$, then (\ref{eq:4}),(\ref{eq:31}) and (\ref{eq:12}) are equivalent
to (\ref{knuutz}).
\end{theorem}

\section{Nonexistence}
\label{sec:Growth-conditions}
We call $u$ a weak super solution of ${\displaystyle -\laplace
  u +g(x,u)=0}$ in $\rz^N$, if $u \in W^{1,2}_{loc}(\rz^N)$ and
\[\int \nabla u \nabla \varphi + \int g(x,u(x))\varphi \ge\, 0 \qquad \forall \phi \in
C^\infty_c(\rz^N),\; \phi \ge 0.\] 
\begin{lemma}
\label{sec:growth-l2}
Suppose (\ref{eq:2}) and 
\begin{align}
\label{eq:52}
\begin{split}
&\text{\ref{eq:1} has a weak positive super solution $u_0$ for some $\lambda>\lambda_1(\rz^N,V)$,}\\
&\text{such that  $u_0 \in C(\rz^N)\cap L^2(\rz^N,V^+) \cap L^p(\rz^N,h^-)$}.    
\end{split}
\end{align}
Then there is a positive solution $u_1 \in
C^1(\rz^N)$ of \ref{eq:1} in $E$ such that $I_\lambda(u_1)<0$ and   
\begin{align}
\label{eq:19}
 (I_\lambda''(u_1)v,v)= \int|\nabla v|^2 - \lambda \int V(x) |v|^2 + (p-1)
\int h(x) |u_1|^{p-2}|v|^2 \ge 0 \quad \forall v \in E.   
\end{align}
\end{lemma}
\begin{proof}
We begin by proving the claim under the additional assumption that 
\begin{align}
\label{eq:22}
u_0 \text{ does not solve the equation in \ref{eq:1}}.  
\end{align}
To this end we consider the set $M:= \left\{u \in E \where 0 \le u(x) \le u_0(x)
\text{ a.e.}\right\}$ and try to minimize $I_\lambda$ in the convex set $M$. Because
$\lambda>\lambda_1(\rz^N,V)$ we may choose $R>0$ large enough such that $\lambda>
\lambda_1(B_R(0),V)$. Because $V \in C(\rz^N)$, a positive minimizer $e_1(B_R(0),V) \in C_c(\rz^N)$
for $\lambda_1(B_R(0),V)$ exists. Due to the fact that $t\cdot e_1(B_R(0),V)$ is an element of $M$
for small values of $t$ and 
$u_0 \in L^2(\rz^N,V^+)\cap L^p(\rz^N,h^-)$ we have $-\infty< \inf_{u \in M}I_\lambda(u)<0$.      
Let $(u_n)_{n \in \nz}$ be a minimizing sequence for $\inf_{u \in M}I_\lambda(u)$. 
Because $0 \le u_n \le u_0$ and $u_0 \in L^2(\rz^N, V^+)\cap L^p(\rz^N,h^-)$ the functions $(u_n)_{n \in \nz}$
are uniformly bounded in $L^2(\rz^N,V^+)\cap L^p(\rz^N,h^-)$.  From the boundedness of
$(I_\lambda(u_n))_{n \in \nz}$ we derive that $(u_n)_{n \in \nz}$ is bounded in $E$. 
We may assume $(u_n)_{n \in \nz}$ converges weakly in $E$ to some $u_1 \in E$ and
pointwise almost everywhere. Since $M$ is convex, $M$ is weakly closed and $u_1 \in M$. 
By the dominated convergence theorem we see that $(u_n)_{n \in \nz}$ converges strongly to $u_1$ in
$L^2(\rz^N,V^+) \cap L^p(\rz^N,h^-)$. Fatou's lemma now yields $I(u_1) \le \inf_{u \in M}
I(u)<0$. By Perron's method (see \cite[Thm. 2.4]{Str96}) the function $u_1$ is a weak solution of $-\laplace
u -V(x) u + h(x) u^{p-1}=0$ in $\rz^N$.\\
By standard regularity results and Harnack's inequality
(see for instance \cite[C.1]{Sim82}) the function $u_1$ is strictly positive and an element of
$C^1(\rz^N)$. (\ref{eq:22}) ensures $u_1 \neq u_0$. Consider the nonnegative function $w:=
u_0-u_1$. Clearly $w$ solves $-\laplace w +\left(- \lambda V(x) + h(x) g(x) \right)w =0
\text{ in }\rz^N,$
where $g$ is a continuous function, defined by
\[ g(x):= \begin{cases}
           \frac{u_1^{p-2}(x)-u_0^{p-2}(x)}{u_1(x)-u_0(x)}& \text{ if }u_1(x)>u_0(x),\\
           (p-2) u_1(x)^{p-3}& \text{ if }u_1(x)=u_0(x).
\end{cases}\]
Hence we may again apply Harnack's inequality to see that $u_1(x)<u_0(x)$ for all $x \in
\rz^N$. Thus $(I_\lambda''(u_1)\varphi,\varphi)\ge 0$ for all $\varphi \in C_c^\infty(\rz^N)$ and by density
for all $v \in E$.\\
Now suppose $u_0$ is a solution of \ref{eq:1}. In this case we consider the equation
\addtocounter{equation}{1} 
\begin{align}
\label{eq:23}
\tag*{$(\theequation)_n$}
-\laplace u - \lambda V(x) u +h(x) u^{p-1}+\frac{1}{n} h^+(x)u^{p-1}=0 \text{ in }\rz^N,
\end{align}
with the corresponding energy functional $J_n$, defined by $J_n(u)=
I_\lambda(u)+\frac{1}{np}\int h^+(x)|u|^p$. The function $u_0$ is a strict super solution of
\ref{eq:23} for each $n \in \nz$. We now apply the same reasoning as above to obtain a sequence of
solutions $(v_n)_{n \in \nz}$ to \ref{eq:23} such that each $v_n$ satisfies (\ref{eq:19})  with $I_\lambda$
replaced by $J_n$. From $J_n(v_n) \le J_1(v_1)<0$ and $0 \le v_n \le u_0$ we conclude that $(v_n)_{n
\in \nz}$ is bounded in $E$. By passing to a subsequence we may assume $(v_n)_{n\in \nz}$ converges
weakly in $E$ to some $u_1 \in E$ and pointwise almost everywhere. By the dominated convergence
theorem and Fatou's lemma we conclude 
\[I_\lambda(u_1) \le \liminf_{n \to \infty}I_\lambda(v_n)\le J_1(v_1)<0.\]
The weak convergence implies that $u_1$ is a critical point of $I_\lambda$. 
Regularity results and Harnack's inequality \cite{Sim82} show $u_1 \in
C^1(\rz^N)$ is a positive solution of \ref{eq:1}.
Fix $w \in E$. By H\"older's
inequality we have for arbitrary $\Omega \subset \rz^N$
\begin{align*}
\int_\Omega |h(x)|\,|v_n|^{p-2}|w|^2 \le \left(\int_\Omega |h(x)|\, |v_n|^p\right)^{\frac{p-2}{p}}
\left(\int_\Omega |h(x)| |w|^p\right)^{\frac{2}{p}} \le C \left(\int_\Omega |h(x)|\,
  |w|^p\right)^{\frac{2}{p}}. 
\end{align*}
We observe that the last term becomes small whenever $|\Omega|$ is small or $\Omega\cap
B_R(0)=\emptyset$ for large $R>0$. Since $(v_n)_{n \in \nz}$ converges pointwise almost
everywhere, we may apply Vitali's theorem 
\cite[Thm. 13.38]{HewStr75} to get 
\[\lim_{n \to \infty}\int h^+(x) |v_n|^{p-2} |w|^2 = \int h^+(x) |u_1|^{p-2} |w|^2 
\text{ and }(I_\lambda''(u_1)w,w) = \lim_{n \to \infty}(J_n''(v_n)w,w) \ge 0.\]   
\end{proof}

\begin{corollary}
Under the assumptions of Lemma \ref{sec:growth-l2} there holds
\begin{align}
\label{eq:33}
D^{1,2}(\rz^N)\cap L^2(\rz^N, V^-) \cap L^p(\rz^N,h^+) \text{ is continuously embedded in }
L^2(\rz^N, V^+).  
\end{align}
\end{corollary}
\begin{proof}
Let $\varphi \in C^\infty_c(\rz^N)$. By (\ref{eq:19}) and H\"older's inequality we have
\begin{align*}
\int V^+(x) |\varphi|^2 &\le \frac{1}{\lambda}\int |\nabla \varphi|^2 + \int V^-(x)|\varphi|^2 + \frac{p-1}{\lambda} \int
h^+(x) u_1^{p-2}|\varphi|^2\\
&\le \frac{1}{\lambda}\| \nabla \varphi\|_2^2 + \|\varphi\|_{L^2(\rz^N,V^-)} +  \frac{p-1}{\lambda} 
\|u_1\|_{L^p(\rz^N,h^+)}^{p-2} \|\varphi\|_{L^p(\rz^N,h^+)}^{2}.     
\end{align*}
The claim follows because $C^\infty_c(\rz^N)$ is dense in the involved spaces.
\end{proof}
The next lemma gives some information concerning the behavior of sequences which converge weakly to zero in
$D^{1,2}(\rz^N)\cap L^p(\rz^N,h^+)$.
\begin{lemma}
\label{sec:growth-l3}
Under the assumptions of Lemma \ref{sec:growth-l2} any sequence $(u_n)_{n \in \nz}$ in $E$, such that 
\begin{align}
\label{eq:25}
u_n \schwach 0 \text{ in } D^{1,2}(\rz^N)\cap L^p(\rz^N,h^+) \text{ and } 
\int |\nabla u_n|^2 \le \frac{\lambda}{2} \int V(x)|u_n|^2,
\end{align}
satisfies $\limsup_{n \to \infty} \int V(x)|u_n|^2 \le 0$.
\end{lemma}
\begin{proof}
Suppose the assertion of the lemma is false. By passing to a subsequence we may assume that      
$(u_n)_{n \in \nz}$ converges to zero pointwise almost everywhere and 
\begin{align}
\label{knutz5}
\lim_{n \to \infty} \int V(x) |u_n|^2 >0.  
\end{align}
By Lemma \ref{sec:growth-l2} there exits $w \in E$ satisfying (\ref{eq:19}). 
We use Vitali's theorem as in the proof of Lemma \ref{sec:growth-l2} to see that  
\begin{align}
\label{eq:27}
\int h^+(x) |w(x)|^{p-2} |u_n|^2 \to 0 \text{ as } n \to \infty.   
\end{align}
By (\ref{eq:25}), (\ref{eq:27}) and (\ref{eq:19}) we have
\begin{align*}
0 \le \int |\nabla u_n|^2 - \lambda \int V(x)|u_n|^2 + (p-1) \int h^+(x) |w|^{p-2} |u_n|^2 \le
-\frac{\lambda}{2} 
\int V(x)|u_n|^2 +o(1), 
\end{align*}
as $n \to \infty$, contrary to (\ref{knutz5}).
\end{proof}
\begin{corollary}
\label{sec:growth-c2}
Under the assumptions of Lemma \ref{sec:growth-l2}, if moreover 
\begin{align}
\label{eq:29}
D^{1,2}(\rz^N) \text{ is continuously embedded in } L^p(\rz^N,h^-),
\end{align}
then we have $\lambda_1(\rz^N,V,h)>0$.  
\end{corollary}
\begin{proof}
Suppose contrary to our claim $\lambda_1(\rz^N,V,h)=0$.
Let $(u_n)_{n \in \nz}$ be a minimizing sequence of $\lambda_1(\rz^N,V,h)$, i.e.
\begin{align}
\label{eq:30}
\int V(x)|u_n|^2 =1,\; \int h(x) |u_n|^p \le 0 \text{ and } \|\nabla u_n\|_2^2 \xpfeil{n \to \infty}
0.  
\end{align}
Clearly $(u_n)_{n \in \nz}$ is bounded in $D^{1,2}(\rz^N)$. We conclude from (\ref{eq:29}) that
$(u_n)_{n \in \nz}$ is bounded in $L^p(\rz^N,|h|)$.
Passing to a subsequence we may assume that $(u_n)_{n \in \nz}$ converges weakly to zero in
$D^{1,2}(\rz^N) \cap L^p(\rz^N,h^+)$. Lemma \ref{sec:growth-l3} tells us that $\limsup_{n \to
  \infty} \int V(x)|u_n|^2 \le 0$ contradicting (\ref{eq:30}).
\end{proof}

\begin{proposition}
\label{sec:growth:p1}
Suppose (\ref{eq:2}) and (\ref{eq:52}). Moreover, we assume  
\begin{align}
\label{eq:8}
D^{1,2}(\rz^N)\cap L^p(\rz^N,|h|) \text{ is compactly embedded in } L^2(\rz^N,V^-),
\end{align}
there exist $\eps <2$, $c_1>0$ and for each $R>2$ there are
 a positive constant $C_R$ and a nonnegative function $h_R$ such that  
\begin{align}
\label{eq:9}
c_1 R^{-\eps} V^+(x) \le V^+(R\cdot x) \text{ for all } x \in \rz^N,\\
\label{eq:10}
h^+(R\cdot x) \le C_R h^+(x) + h_R(x)\text{ for all } x \in \rz^N \text{ and}\\
\label{eq:11}
D^{1,2}(\rz^N) \cap L^p(\rz^N,h^+) \text{ is continuously embedded in } L^p(\rz^N, h_R).  
\end{align}
Then $D^{1,2}(\rz^N) \cap L^p(\rz^N,h^+)$ is compactly embedded in
$L^2(\rz^N,V^+)$. 
\end{proposition}
\begin{proof}
For $u \in D^{1,2}(\rz^N) \cap L^p(\rz^N,h^+)\cap L^2(\rz^N,V^+)$ we write $u_R(x):= R^{(\eps-N)/2}
u(x/R)$ and deduce by (\ref{eq:9}) and (\ref{eq:10}) for $R>2$
\begin{align}
\label{eq:35}
\int V^+(x) |u|^2 \le c_1^{-1} R^\eps \int V^+(R\cdot x) |u|^2 = c_1^{-1}\int V^+(x) |u_R|^2,\\
\label{eq:36}
\int h^+(x) |u_R|^p \le R^{\frac{p\eps}{2}-N \frac{p-2}{2}} \int h^+(R \cdot x) |u|^p \le C'_R \int
h^+(x) |u|^p + C'_R \int h_R(x) |u|^p,\\
\label{eq:37}
\int |\nabla u_R|^2 = R^{\eps -2} \int |\nabla u|^2. 
\end{align}Let $(u_n)_{n \in \nz}$ be a sequence in $D^{1,2}(\rz^N) \cap L^p(\rz^N,h^+)$ converging weakly to
zero. The assertion of the theorem follows if we show that $(u_n)_{n \in \nz}$ converges to zero in
$L^2(\rz^N,V^+)$. Conversely, by passing to a subsequence, suppose that
\begin{align}
\label{eq:38}
\lim_{n \to \infty}\int V^+(x) |u_n|^2 \text{ exists and is positive.}  
\end{align}
First we show that (\ref{eq:38}) leads to a contradiction under the additional assumption that
\begin{align}
\label{eq:39}
\int |\nabla u_n|^2 \le \frac{\lambda}{2}\int V^+(x)|u_n|^2 \text{ for all }n \in \nz.    
\end{align}
By (\ref{eq:8}) the sequence of integrals $\int V^-(x) |u_n|^2$ converges to zero. Consequently Lemma
\ref{sec:growth-l3} and (\ref{eq:39}) yield $\int V^+(x) |u_n|^2 \to 0$ as $n \to \infty$, contrary
to (\ref{eq:38}).\\
To complete the proof we suppose 
\begin{align*}
\limsup_{n \to \infty}\frac{\int |\nabla u_n|^2}{\int V^+(x) |u_n|^2}=: \Lambda\ge \frac{\lambda}{2}.
\end{align*}
If $\Lambda = \infty$, then obviously $\int V^+(x) |u_n|^2 \to 0$ as $n \to \infty$ because
$\|\nabla u_n\|_2$ is uniformly bounded in $n$.
If $\Lambda <\infty$, then we take $R\ge (3\Lambda/\lambda c_1)^{1/(2-\eps)}$ and  
\begin{align}
\label{eq:32}
v_n(x):= u_{n,R}(x)=R^{(\eps-N)/2} u_n(x/R) \text{ for each } n \in \nz.  
\end{align}
By (\ref{eq:36}) and (\ref{eq:37}) the sequence $(v_n)_{n \in \nz}$ is bounded in $D^{1,2}(\rz^N)\cap
L^p(\rz^N,h^+)$ and converges weakly to zero in $D^{1,2}(\rz^N)$. Thus by passing to a subsequence we
may assume $(v_n)_{n \in \nz}$ converges weakly to zero in $D^{1,2}(\rz^N)\cap L^p(\rz^N,h^+)$. From
(\ref{eq:35}) and (\ref{eq:37}) we deduce
\begin{align}
\label{eq:53}
\limsup_{n \to \infty}\frac{\int |\nabla v_n|^2}{\int V^+(x) |v_n|^2} \le \frac{\lambda c_1}{3
  \Lambda c_1} \frac{\int |\nabla u_n|^2}{\int V^+(x) |u_n|^2} \le \frac{\lambda}{3}.  
\end{align}
From the above analysis we conclude $\int V^+(x) |v_n|^2 \to 0$ as $n \to \infty$, hence by
(\ref{eq:35}) we have $\int V^+(x) |u_n|^2 \to 0$ as $n \to \infty$,
contrary to (\ref{eq:38}).
\end{proof}
\begin{remark}
\label{sec:growth:rem1}
Using (\ref{eq:32}) and (\ref{eq:53}) in the proof of Proposition \ref{sec:growth:p1}
we easily obtain that $\lambda_1(\rz^N,V)=0$ under the assumption (\ref{eq:9}). 
\end{remark}
\begin{corollary}
\label{sec:growth-conditions-c3}
Suppose (\ref{eq:2}),(\ref{eq:52}) and (\ref{eq:28})-(\ref{eq:40}).
Then $D^{1,2}(\rz^N) \cap L^p(\rz^N,h^+)$ is compactly embedded in $L^2(\rz^N,V^+)$.  
\end{corollary}
\begin{proof}
We may estimate for $R \ge 1$ by (\ref{eq:40})
\begin{align}
\label{eq:50}
h^+(R \cdot x) &\le c_R h^+(x) + \|h^+ \cdot\charak_{B_{M\cdot R}(0)} \|_\infty \cdot
\charak_{B_{M\cdot R}(0)}(x) \text{ for some }M>0.
\end{align}
By (\ref{eq:28}) we may find positive constants $C_2>2$, $0<c_1^-<c_1^+$ such that
\begin{align}
\label{eq:49}
c_1^-|x|^\alpha \le V(x) \le c_1^+ |x|^\alpha \quad \text{ for all } |x|\ge C_2.    
\end{align}
Define $V_0 :\rz^N \to \rz$ by
\begin{align*}
V_0(x) := \begin{cases}
               c_1^-|x|^\alpha &\text{ if } C_2+1\le |x|\\
               c_1^-(|x|^\alpha - C_2^\alpha) &\text{ if } C_2\le |x|\le C_2+1\\
               \min(V(x),0) &\text{ if } |x|\le C_2.
               \end{cases}
\end{align*}
We may estimate for all $|x|\ge C_2+1$ and $R\ge 1$
\begin{align*}
\frac{V_0(x)}{V_0(Rx)} &\le R^{-\alpha}. 
\end{align*}
For $C_2 \le |x| \le C_2+1$ and $R\ge 2$ we have since $C_2>2$
\begin{align*}
\frac{V_0(x)}{V_0(Rx)} &\le \frac{c_1^-(|C_2+1|^\alpha - C_2^\alpha)}{c_1^- R^\alpha C_2^\alpha}
\end{align*}
Putting together the two above observations we obtain for all $R>2$
\begin{align}
\label{eq:51}
\frac{C_2^\alpha}{(|C_2+1|^\alpha - C_2^\alpha)} R^\alpha {V_0}^+(x)\le {V_0}^+(Rx)    
\end{align}
By Remark \ref{sec:growth:rem1} we obtain $\lambda_1(\rz^N,V_0)=0$. Since $V_0 \le V$ 
and (\ref{eq:50})-(\ref{eq:51}) hold, Proposition \ref{sec:growth:p1} applied to $V_0$ and $h$
gives
\begin{align}
\label{eq:54}
D^{1,2}(\rz^N) \cap  L^p(\rz^N,h^+) \text{ is compactly embedded in } L^2(\rz^N,V_0^+).    
\end{align}
Since $V^+$ is smaller than $(c_1^-)^{-1}c_1^+V_0^+$ outside a bounded domain, (\ref{eq:54})
remains true if we replace $V_0^+$ by $V^+$. 
\end{proof}
To prove that (\ref{eq:31})-(\ref{eq:12}) are equivalent to (\ref{knuutz}) we need
\begin{lemma}
\label{sec:growth:l1}
Under the assumptions of Lemma \ref{sec:growth-l2}, if moreover there exits $R>0$ such that 
\begin{align}
\label{eq:24}
\begin{split}
V,\,h \in C^2(\rz^N \backslash \overline{B_R(0)}),\quad h(x) >0 \; \forall |x|\ge R,\quad
\limsup_{|x| \to \infty} \frac{V^+(x)}{h(x)} = 0,\\
\laplace \left(\frac{V^+(x)}{h(x)}\right)^{\frac{1}{p-2}}\ge 0\;
\text{ in } \Omega_+ := \{x \in \rz^N \where V(x)>0,\, |x|>R\}.
\end{split}
\end{align}
Then the integral ${\displaystyle \int_{\rz^N \backslash B_R(0)} V^+
  \left(\frac{V^+}{h}\right)^{\frac{2}{p-2}}}$ is finite.
\end{lemma}
\begin{proof}
By Lemma \ref{sec:growth-l2} and standard regularity results we may assume $u_0$ solves \ref{eq:1}
and $u_0 \in C^2(\rz^N \backslash \overline{B_R(0)}) \cap E$. For abbreviation we write
$w(x)$ instead of $(V^+(x)/h(x))^{1/(p-2)}$. Because $u_0$ is a positive 
function we may choose $\sigma >0$ such that $\lambda>\sigma$ and 
\[u_0(x)> (\lambda - \sigma)^{\frac{1}{p-2}} w(x) \text{ in } \{x \in \rz^N \where |x|=R\}.\] 
We will denote by $\Omega_\sigma$ the set $\{x \in \rz^N \where -\laplace u_0
-\sigma V(x) u_0 >0,\, |x|>R\}$. Because $u_0$ is a solution to \ref{eq:1} an easy calculation shows
\begin{align*}
x \in \Omega_\sigma \aequi h(x) u_0(x)^{p-2}< (\lambda -\sigma) V(x) \aequi
u_0(x)<(\lambda-\sigma)^{\frac{1}{p-2}} w(x)  \text{ and }
\overline{\Omega_\sigma} \subset \Omega_+.  
\end{align*}
The basic idea of the proof is to show that $\Omega_\sigma$ is empty.  Conversely, suppose that
there is $x_0 \in \Omega_\sigma$, i.e. $|x_0|>R$ and 
\begin{align}
\label{eq:41}
u_0(x_0) -(\lambda-\sigma)^{1/(p-2)} w(x) = -\eps \text{ for some } \eps>0.   
\end{align}
From (\ref{eq:24}) we may choose $R_1>|x_0|$ large enough to have $u_0(x)
-(\lambda-\sigma)^{1/(p-2)} w(x)>-\eps/2$ for all $x \in \overline{\Omega_\sigma}$ satisfying
$|x|=R_1$. With the notation $\Omega_* := \Omega_\sigma \cap B_{R_1}(0)$ we have
\begin{align}
\label{eq:15}
 \inf_{x \in \rand \Omega_*}\left(u_0(x) -(\lambda-\sigma)^{1/(p-2)} w(x)\right) \ge -\frac{\eps}{2}
\text{ and } \laplace(u_0-(\lambda-\sigma)^{1/(p-2)} w) \le 0 \text{ in } \Omega_*.  
\end{align}
By (\ref{eq:15}) and the maximum principle we deduce $u_0(x) -(\lambda-\sigma)^{1/(p-2)} w(x) \ge
-\eps/2$ for all $x \in \Omega_*$, contrary to (\ref{eq:41}).\\
Consequently, we have $u_0(x) \ge (\lambda-\sigma)^{1/(p-2)} w(x)$ for all $x$ satisfying $|x| \ge R$
and we may estimate
\[\infty > \int_{\rz^N \backslash B_R(0)} V^+(x)u_0^2 \ge (\lambda -\sigma)^{\frac{2}{p-2}}   
\int_{\rz^N \backslash B_R(0)} V^+\left(\frac{V^+}{h}\right)^{\frac{2}{p-2}}.\]
\end{proof}

\begin{corollary}
\label{sec:growth-c1}
Under the hypothesis of Lemma \ref{sec:growth-l2}, given any functions $h_0$ and $V_0$ satisfying (\ref{eq:24}) such that
\begin{align}
\label{eq:42}
h_0(x) \ge h(x) \text{ and } V_0(x) \le V(x) \text{ in } \rz^N,\\
\label{eq:43}
\exists 0 \lneqq \varphi \in C^\infty_c(\rz^N): \int |\nabla \varphi|^2 - \lambda \int V_0(x) |\varphi|^2 < 0,
\end{align}
then the integral ${\displaystyle \int_{\rz^N \backslash B_R(0)} V_0^+ \left(\frac{V_0^+}{h_0}\right)^{\frac{2}{p-2}}}$
is finite for some $R>0$. 
\end{corollary}
\begin{proof}
By (\ref{eq:42}) the function $u_0$ is a super solution of $-\laplace u -V_0(x) u + h_0(x) u^{p-1}=0$ in $\rz^N$. 
From (\ref{eq:43}) we conclude $\lambda > \lambda_1(\rz^N,V_0)$. Hence we may apply Lemma
\ref{sec:growth:l1} with $h$ and $V$ replaced by $h_0$ and $V_0$.
\end{proof}
\begin{proposition}
\label{sec:growth:p2}
Suppose (\ref{eq:2}),(\ref{eq:32}), (\ref{eq:28})-(\ref{eq:40}) and $(N-4)\le \alpha$, where
$\alpha$ is given in (\ref{eq:28}). 
Then we have for some $R>0$
\[\int_{\rz^N \backslash B_R(0)} V^+
  \left(\frac{V^+}{h}\right)^{\frac{2}{p-2}}< \infty.\]
\end{proposition}
\begin{proof}
An easy computation shows that the claim of the proposition holds if and only if 
$(\beta-\alpha)>(N/2+\alpha/2)(p-2)$. Conversely, suppose that $(\beta-\alpha)\le
(N/2+\alpha/2)(p-2)$. Since $(N-4)\le \alpha$ there holds 
\begin{align}
\label{eq:16}
(N-2)(p-2) \le \left(\frac{N}{2}+\frac{\alpha}{2}\right)(p-2). 
\end{align}
By direct calculation we obtain
\begin{align}
\label{eq:17}
\laplace\left(|x|^{-\frac{\beta-\alpha}{p-2}}\right) \ge 0 \text{ for all }x \neq 0 \text{ if } 
(\beta-\alpha) \ge (N-2)(p-2).  
\end{align}
From (\ref{eq:16}) it is possible to choose $\beta' \ge \beta$ such that 
\begin{align}
\label{eq:18}
(N-2)(p-2) \le (\beta'-\alpha) \le (N/2+\alpha/2)(p-2).   
\end{align}
Consider the functions $h_0$ and $V_0$ defined by
\begin{align*}
h_0(x) := \max(C_2 |x|^{\beta'}\charak_{|x|>R}+C_2 R^{\beta'}\charak_{|x|\le R},h(x)),\\
V_0(x) := \min(C_1 |x|^{\alpha}\charak_{|x|>R}+C_1 R^{\alpha}\charak_{|x|\le R},V(x)).  
\end{align*}
Take an arbitrary nonnegative $\varphi \in C_c^\infty(\rz^N \backslash \overline{B_R(0)})$ and
define for $t>0$ the function $\varphi_t(x) := t^{(2-N)/2} \varphi(x/t)$. 
Then $\|\nabla \varphi_t\|_2 = \|\nabla \varphi\|_2$ for all $t$ and since $\alpha >-2$
\begin{align}
\label{eq:58}
\int_{\rz^N \backslash B_R(0)} |x|^\alpha |\varphi_t(x)|^2 =
t^{2+\alpha} \int_{\rz^N \backslash B_{R/t}(0)} |x|^\alpha |\varphi(x)|^2 \to \infty \text{ as
  } t \to \infty.   
\end{align}
Hence $\lambda_1(\rz^N,V_0)=0$ and we may apply Corollary \ref{sec:growth-c1} to see that
\[\int_{\rz^N \backslash B_R(0)} V_0^+ \left(\frac{V_0^+}{h_0}\right)^{\frac{2}{p-2}}<\infty.\]
Thus $(\beta'-\alpha) > (N/2+\alpha/2)(p-2)$, contrary to (\ref{eq:18}).
\end{proof}
\begin{proof}[\bf Proof of Theorem \ref{sec:intro:t8}]
Analysis similar to that in the proof of Proposition \ref{sec:growth:p2} in (\ref{eq:58}) leads to
$\lambda_1(\rz^N,V)=0$.\\
Obviously (\ref{eq:12}) is a consequence of (\ref{eq:31}).
If (\ref{eq:12}) holds, then Corollary \ref{sec:growth-conditions-c3}
shows that (\ref{eq:4}) is satisfied. If 
(\ref{eq:4}) holds, then by Lemma \ref{s>0l1} below we have $\lambda_1(\rz^N,V,h)>0$. Consequently
Theorem \ref{sec:intro:t1} yields (\ref{eq:31}).
The last part of our claim is true by applying Proposition \ref{sec:growth:p2} and since
(\ref{knuutz}) is a sufficient condition for (\ref{eq:4}).      
\end{proof}
\section*{Spectral conditions} 
The contents of the next lemma is the relation of positive solutions to the first eigenvalue of
linear problems. It is known and can be derived by testing 
${-\laplace u_0 = k(x) u_0}$ with ${\varphi^2/u_0}$ and 
Cauchy-Schwarz' inequality. 
\begin{lemma}  
\label{snonl1}  
Suppose there exists a continuous, positive weak super solution $u_0>0$ of 
${\displaystyle -\laplace u = k(x) u}$ in
$\rz^{N}$, where ${\displaystyle k \in L^1_{loc}(\rz^N)}$. Then 
\begin{equation}  
\int k(x) \varphi^{2} \le \int |\nabla \varphi|^{2} \quad \forall \; \varphi \in
C^\infty_c(\rz^{N}).  
\label{snone1}  
\end{equation}  
\end{lemma}  
\begin{proof}[\bf{Proof of Theorem \ref{sec:intro:t5}}]
Suppose $(\lambda,u)$ solves \ref{eq:1}. The fact, that $C_c^\infty(\rz^N)$ is dense in the space
$E_u:=D^{1,2}(\rz^N)\cap L^2(\rz^N,|V|)\cap L^2(\rz^N,h^+ u^{p-2})$ and  Lemma \ref{snonl1} imply
\begin{align}
\label{eq:20}
\begin{split}
\frac{1}{\lambda} &\ge \sup_{\substack{\phi \in C_c^\infty(\rz^N)\\\phi \neq 0}} 
\frac{\int \left(V(x) - \frac{h(x) u^{p-2}}{\lambda}\right)\phi^2}{\|\nabla \phi\|_2^2}
\ge \sup_{\substack{\phi \in C_c^\infty(\rz^N)\\\phi \neq 0}} 
\frac{\int \left(V(x) - \frac{h^+(x) u^{p-2}}{\lambda}\right)\phi^2}{\|\nabla \phi\|_2^2}\\
&= \sup_{\substack{\phi \in E_u  \\ \phi \neq 0}} 
\frac{\int \left(V(x) - \frac{h^+(x) u^{p-2}}{\lambda}\right)\phi^2}{\|\nabla \phi\|_2^2}
\ge \sup_{\substack{\phi \in D^{1,2}(\rz^N)\cap L^2(\rz^N,|V|)\\\phi \neq 0\\
\phi(x)=0 \text{ a.e. in }\rz^N \backslash \Omega^{-0}}} 
\frac{\int V(x)\, \phi^2}{\|\nabla \phi\|_2^2}\\ 
&= \frac{1}{\lambda_1(\Omega^{-0},V)}.
\end{split}   
\end{align}
Now suppose $\lambda_1(\Omega^{-0},V)<\infty$ is attained by $v_0$ and $u$ is a solution of
\ref{eq:1} for $\lambda= \lambda_1(\Omega^{-0},V)$. The above
calculation shows for $\lambda=\lambda_1(\Omega^{-0},V)$
\begin{align}
\label{eq:21}
\frac{1}{\lambda_1(\Omega^{-0},V)}= 
\sup_{\substack{\phi \in D^{1,2}\cap L^2(|V|)\cap L^2(h^+ u^{p-2})  \\ \phi \neq 0}} 
\frac{\int \left(V(x) - \frac{h^+(x) u^{p-2}}{\lambda_1(\Omega^{-0},V)}\right)\phi^2}{\|\nabla \phi\|_2^2}
= \frac{\int V(x) \,v_0^2}{\|\nabla v_0\|_2^2}.  
\end{align}
Consequently $v_0$ weakly solves 
\[-\laplace v_0(x) -\lambda_1(\Omega^{-0},V) V(x) v_0(x) + h^+(x) u^{p-2}(x)v_0(x)=0 \qquad x \in \rz^N.\]
Harnack's inequality and a regularity result in \cite[C.1]{Sim82} imply $v_0(x)>0$ for all $x \in \rz^N$. 
Hence $h^+ \equiv 0$ and $\Omega^{-0}=\rz^N$.\\
Since $C_c^\infty(\rz^N)$ is dense in $D^{1,2}(\rz^N)\cap
L^2(\rz^N,|V|)$ and by (\ref{eq:20})-(\ref{eq:21}) 
\[\frac{1}{\lambda_1(\Omega^{-0},V)}= \sup_{\substack{\phi \in C_c^\infty(\rz^N)\\\phi \neq 0}} 
\frac{\int \left(V(x) - \frac{h(x) u^{p-2}}{\lambda_1(\Omega^{-0},V)}\right)\phi^2}{\|\nabla \phi\|_2^2}=
\sup_{\substack{\phi \in C_c^\infty(\rz^N)\\\phi \neq 0}} 
\frac{\int V(x) \phi^2}{\|\nabla \phi\|_2^2} =\frac{1}{\lambda_1(\Omega^{-0},V)}.\]
Hence $\int h^-(x) u^{p-2} v_0^2 =0$ and $h^- \equiv 0$, which is impossible. 
\end{proof}

\section{Behavior of ${\bf \lambda_1}$}
\label{sec:behav}
\begin{lemma}
\label{s>0l1}
Suppose (\ref{eq:2})-(\ref{eq:4}) hold and $\lambda_1(\rz^{N},V,h)$ is finite.
Then $\lambda_1(\rz^{N},V,h)$ is attained by some nonnegative $u_0 \in E$.
\end{lemma}
\begin{proof}
Let $(u_n)_{n \in \nz}$ be a minimizing sequence for $\lambda_1(\rz^{N},V,h)$ . Because the inclusions of
$D^{1,2}(\rz^{N})$ in $L^{p}(\rz^N,h^{-})$ and $D^{1,2}(\rz^{N})\cap L^p(\rz^N,h^+)$ in $L^2(\rz^N,V^+)$
are compact, we see that $(u_n)_{n \in \nz}$ is bounded in $E$. 
Hence we may assume $(u_n)_{n \in \nz}$ converges weakly in $E$ to some
$u_0\in E$. Due to the compactness of the above embeddings we have
$\int V(x)|u_0|^2 \ge 1$, $\int h(x) |u_0|^{p}\le 0$ and $\|\nabla u_0\|^{2} \le \lambda_1(\rz^N,V,h)$.
The rest of the claim follows after appropriate scaling and replacing $u_0$ by $|u_0|$.
\end{proof}

\begin{proof}[\bf{Proof of Theorem \ref{sec:intro:l1}}]
By definition $\lambda_1(\rz^{N},V,h_\mu) \le \lambda_1(\Omega^{-0},V)$ for all $\mu\ge 1$ and the value
$\lambda_1(\rz^{N},V,h_\mu)$ is monotone increasing in $\mu$. Therefore
\[ \lim_{\mu \to \infty} \lambda_1(\rz^{N},V,h_\mu) =: \Lambda \le \lambda_1(\Omega^{-0},V).\]
If $\Lambda=+\infty$ then obviously the claim is true. Hence we may assume $\Lambda <\infty$.
Consider ${\displaystyle (u_n)_{n \in \nz}}$, where
$u_{n}$ is the minimizer of $\lambda_1(\rz^{N},V,h_n)$ found in Lemma \ref{s>0l1}. Analysis similar to
that in the proof of Lemma \ref{s>0l1} yields that $(u_{n})_{n \in \nz}$ is bounded in $E$ and by passing
to subsequence we may assume
\[ u_n \schwach u_0 \text{ in } E, \quad \int V^+(x)|u_n|^2 \to \int V^+(x)|u_0|^2 \text{ and }
\int V |u_0|^2 \ge 1. \]
Because ${\int h^{-}(x)|u_n|^{p}}$ is bounded, we see
${\int h^{+}(x) |u_n|^{p} \le \frac{C}{n}}$. From that we conclude
\[0 \le \int h^{+}(x) |u_0|^{p} \le \liminf_{n \to \infty} \int h^{+}(x)|u_n|^{p} \le 0,\]
hence that $u_0(x)=0$ for almost every $x \in \rz^{N} \backslash \Omega^{-0}$, and finally that
\[ \lambda_1(\Omega^{-0},V) \le \|\nabla u_0\|_2^{2} \le \liminf_{n \to \infty} \|\nabla
u_n\|_2^{2} = \lim_{n \to \infty} \lambda_1(\rz^N,V,h_{n}) = \Lambda \le \lambda_1(\Omega^{-0},V).\]
\end{proof}

\section{Existence of a local minimizer}
\label{sec:Exist-local-minim}
The existence of a minimizer in an open subset of $E$ is proved by
\begin{lemma}
\label{s>0l4}
Suppose (\ref{eq:2})-(\ref{eq:4}) and $\lambda_1(\rz^N,V)< \lambda < \lambda_1(\rz^N,V,h)$
hold. Then the value
\begin{align}
\label{eq:26}
\sigma(\lambda) := \inf\left\{I_{\lambda}(u)\where \frac{\int V(x) |u|^2}{\|\nabla u\|_2^{2}}
> \frac{1}{\lambda}, \; u \in E \backslash \{0\}\right\}
\end{align}
is negative and attained by some nonnegative $u_0 \in E$.
\end{lemma}
\begin{proof}
Because $\lambda_1(\rz^N,V)<\lambda$, there is a $\tilde u \in E$ such that
\[\|\nabla \tilde u\|_2^{2}-\lambda \int V |\tilde u|^2 <0.\]
Hence we have $I_{\lambda}(t\tilde u)<0$ for small positive $t$ and $\sigma(\lambda)
<0$.\\
Let $(u_n)_{n \in \nz}$ be a minimizing sequence. To obtain a contradiction, suppose
\begin{align}
\label{eq:55}
0\le s_n^2:=\int V(x) |u_n|^2 \to \infty \text{ as } n \to \infty.  
\end{align}
Write ${v_n:= \frac{u_n}{s_n}}$. The sequence $(v_n)_{n \in \nz}$ is bounded in
$D^{1,2}(\rz^N)$ by (\ref{eq:26}) and we conclude for large $n$
\begin{align*}
0 &>\frac{1}{2}\|\nabla v_n\|^{2} - \frac{\lambda}{2} + \frac{s_n^{p-2}}{p}\int h(x) |v_n|^{p}
\end{align*}
Consequently we have $\limsup_{n \to \infty} \int h(x) |v_n|^{p} \le 0$. 
By (\ref{eq:3})-(\ref{eq:4}) we see $(v_n)_{n \in \nz}$ is bounded in $E$ and by passing to a
subsequence and (\ref{eq:3}) we may assume 
\begin{align}
\label{eq:56}
\int h(x) |v_0|^{p} \le \liminf_{n \to \infty} \int h(x) |v_n|^{p} \le 0.  
\end{align}
(\ref{eq:4}) leads to $\int V(x) |v_0|^2 \ge 1$.
Finally ${\|\nabla v_0\|_2^{2} \le \liminf_{n \to \infty} \|\nabla v_n\|_2^{2} \le
\lambda < \lambda_1(\rz^N,V,h)}$ contradicts the definition of $\lambda_1(\rz^N,V,h)$.\\
Hence $(u_n)_{n \in \nz}$ is bounded in $D^{1,2}(\rz^N)$, therefore $ \int h^{-}(x)
|u_n|^{p}$ is bounded due to (\ref{eq:3}). Because
$I_{\lambda}(u_n)<0$ for large $n$ we may deduce that 
$\int h^{+}(x) |u_n|^{p}$ is bounded. From (\ref{eq:4}) we infer
that $\int |V(x)| |u|^2$ is bounded, consequently $(u_n)_{n \in \nz}$ is
bounded in $E$. Taking a subsequence we may assume $(u_n)_{n \in \nz}$ converges weakly in $E$
and strongly in  $L^2(\rz^N,V^+)$ and $L^p(\rz^N,h^-)$ to some $u_0
\in E$. Thus 
\[I_{\lambda}(u_0) \le \liminf_{n \to \infty} I_{\lambda}(u_n)=
\sigma(\lambda)<0 \text{ and }  \|\nabla u_0\|_2^{2}- \lambda \int
V(x) |u_0|^2\le 0.\]
If $\|\nabla u_0\|_2^{2}=\lambda \int V(x) |u_0|^2$ then due to the
definition of $\lambda_1(\rz^N,V,h)$ we have $\int h(x) |u_0|^p \ge 0$
and consequently $I_{\lambda}(u_0)\ge 0$ which is impossible.
Hence $|u_0|$ minimizes $\sigma(\lambda)$.
\end{proof}

\begin{proof}[\bf{Proof of Theorem \ref{sec:intro:t1}}]
The part of our claim concerning existence is an obvious consequence
of Lemma \ref{s>0l4}.\\
Suppose contrary to our claim that $h\ge 0$ and $u \neq v$ are two positive
solutions of \ref{eq:1} in $C^1(\rz^N)\cap E$. A direct 
computation shows that $f\in C(\rz^N)$, defined by
\[f(x):= 
\begin{cases} 
\frac{u^{p-1}(x)-v^{p-1}(x)}{u(x)-v(x)} &\text{if } u(x)\neq v(x),\\
(p-1)u^{p-2}(x) &\text{otherwise.}
\end{cases}\]
satisfies $f(x)>u^{p-2}(x)$ for all $x \in \rz^N$. We write $w= u-v$.
Suppose $w(x_0)\neq 0$ for some $x_0 \in \rz^N\backslash \Omega^{-0}$. Then we obtain a contradiction
by using (\ref{snone1}) and
\begin{align*}
0 &= \left(I'_\lambda(u)-I'_\lambda(v)\right)w = \int |\nabla w|^2 - \lambda \int V(x) w^2 + \int
h(x) f(x) w^2\\
&> \int |\nabla w|^2 - \lambda \int V(x) w^2 + \int h(x) u^{p-2} w^2 \ge 0.   
\end{align*}
Hence $w(x)=0$ for all $x \in \rz^N \backslash \Omega^{-0}$ and $\|\nabla w\|_2^2=\lambda
\int V(x) w^2$, contrary to $\lambda< \lambda_1(\Omega^{-0},V)$.
\end{proof}

\begin{remark}
Under the assumptions of Theorem \ref{sec:intro:t1} we obtain a solution $u_\lambda$ of \ref{eq:1}
for every $\lambda_1(\rz^N,V)<\lambda<\lambda_1(\rz^N,V,h)$ as local minimizer of $I_\lambda$, i.e.
\[I_\lambda(u_\lambda) = \inf\left\{ I_\lambda(u) \,\where\, \frac{\int V(x) |u|^2}{\|\nabla
u\|_2^{2}} > \frac{1}{\lambda}, \; u \in E\backslash\{0\}\right\}.\]
Consequently we have
\begin{align}
\label{eq:176}
I_{\lambda'}(u_{\lambda'}) \le I_{\lambda}(u_{\lambda})-(\lambda'-\lambda) \int V(x) |u_\lambda|^2
\le I_\lambda(u_\lambda) \text{ if } \lambda \le \lambda'.   
\end{align}
\end{remark}

\begin{proof}[\bf{Proof of Theorem \ref{sec:intro:t2}}]
To obtain a contradiction we suppose $\lambda_1(\rz^N,V)=\lambda_1(\Omega^{-0},V)$. Obviously
$\lambda_1(\rz^N,V)$ is finite. Hence by Lemma \ref{sec:intro:l1} the value
$\lambda_1(\Omega^{-0},V)$ is attained by some  
function $e_1(\Omega^{-0},V)$. Because $\lambda_1(\rz^N,V)=\lambda_1(\Omega^{-0},V)$ the function
$e_1(\Omega^{-0},V)$ is a weak solution of $-\laplace u= \lambda_1(\rz^N,V) V(x) u$ in
$\rz^N$. Harnack's inequality yields $e_1(\Omega^{-0},V)$ is positive in $\rz^N$, which is
impossible. Existence and uniqueness follow from Theorem \ref{sec:intro:t1}.\\
We observe for any $\mu < \lambda_1(\Omega^{-0},V)$ 
\begin{align}
\label{eq:171}
\sup_{\lambda_1(\rz^N,V)<\lambda\le \mu} \|u_\lambda\|_E <\infty.  
\end{align}
To show (\ref{eq:171}) we proceed analogously to the proof of Lemma \ref{s>0l4} and assume,
contrary to (\ref{eq:171}), that there is a sequence $(\lambda_n)_{n \in \nz}$ such that 
$\|u_{\lambda_n}\|_E \to \infty$ and $\lambda_n \to \mu <\lambda_1(\Omega^{-0},V)$ 
as $n \to \infty$. Because $h$ is nonnegative and $I_{\lambda_n}(u_{\lambda_n})<0$ we obtain 
$s_n^2:= \int V(x) |u_{\lambda_n}|^2 \to \infty$ as $n \to \infty$. We write $v_n :=
u_{\lambda_n}/s_n$ and get
\begin{align}
\label{eq:172}
0 \ge \frac{I_{\lambda_n}(u_{\lambda_n})}{s_n^2} = \frac{1}{2}\| v_n\|^2 - \frac{\lambda_n}{2} +
\frac{s_n^{p-2}}{p} \int h(x) |v_n|^p. 
\end{align}
Consequently $(v_n)_{n \in \nz}$ is bounded in $E$ and passing to a subsequence we may assume
$(v_n)_{n \in \nz}$ converges weakly in $E$ to $v_0$. We have $v_0 \neq 0$, because 
\[\int V(x)|v_0|^2 \ge \limsup_{n \to \infty} V(x)|v_n|^2 =1.\]
Finally $\int h(x)|v_0|^p=0$ and $\|\nabla v_0\|_2^2 \le \mu$ contradict
the definition of $\lambda_1(\Omega^{-0},V)$.\\
Fix $(\lambda_n)_{n \in \nz}$ such that $\lambda_n \to \lambda \in
\left(\lambda_1(\rz^N,V),\lambda_1(\Omega^{-0},V)\right)$ as $n \to \infty$. To obtain a
contradiction we suppose 
\begin{align}
\label{eq:230}
\|u_{\lambda_n}-u_\lambda\|_E \ge \eps \; \text{ for all }  n \in \nz \text{ and some }\eps >0.  
\end{align}
From (\ref{eq:176}) we infer
\begin{align}
\label{eq:231}
\limsup_{n \to \infty} I_{\lambda_n}(u_{\lambda_n}) < 0.  
\end{align}
By (\ref{eq:171}) the sequence $(u_{\lambda_n})_{n \in \nz}$ is bounded in $E$. Passing to a
subsequence we may assume $(u_{\lambda_n})_{n \in \nz}$ converges weakly in $E$ and strongly in
$L^2(\rz^N,V^+)$ to some $u_0 \in E$. The function $u_0$ is a weak nonnegative solution of
\ref{eq:1} and $u_0 \not\equiv 0$ because
\[I_\lambda(u_0) \le \liminf_{n \to \infty} I_{\lambda_n}(u_{\lambda_n})<0.\]
Uniqueness of positive solution gives $u_0=u_\lambda$. We may estimate using the strong convergence in
$L^2(\rz^N,V^+)$,
\begin{align*}
0 &= I_{\lambda_n}'(u_{\lambda_n})u_{\lambda_n} -I_{\lambda}'(u_{\lambda})u_{\lambda}\\
&= \left(\|\nabla u_{\lambda_n}\|_2^2-\|\nabla u_{\lambda}\|_2^2\right) + 
\lambda_n\left(\int V^-(x) |u_{\lambda_n}|^2\!-\!\int V^-(x)
  |u_{\lambda}|^2\right) + 
\underbrace{(\lambda_n\!-\!\lambda) \int V^-(x) |u_\lambda|^2}_{\to 0 \text{ as }n \to \infty}\\
&\quad - \underbrace{\left(\lambda_n
  \int V^+(x)|u_{\lambda_n}|^2 - \lambda \int V^+(x)|u_{\lambda}|^2\right)}_{\to 0 \text{ as }n \to \infty}
 + \left(\int h(x)|u_{\lambda_n}|^p -\int h(x)|u_{\lambda}|^p\right).       
\end{align*}
Consequently, the uniform convexity of the involved spaces shows that $u_{\lambda_n}\to u_\lambda$ in
$E$ as $n \to \infty$, contradicting (\ref{eq:230}).\\ 
Fix $(\lambda_n)_{n \in \nz}$ such that $\lambda_n \to
\lambda_1(\rz^N,V)$ as $n \to \infty$. By (\ref{eq:171}) we see 
\begin{align}
\label{eq:175}
0 \ge I_\lambda(u_\lambda)= \frac{1}{2}\|\nabla u_\lambda \|_2^2- \lambda_1(\rz^N,V) \int
V(x) |u_\lambda|^2 + \frac{1}{p} \int h(x)
|u_\lambda|^p + o(1)_{\lambda \to \lambda_1(\rz^N,V)}.  
\end{align}
Hence $u_{\lambda_n} \to 0$ in $L^p(\rz^N,h)$ as $n \to \infty$. 
Suppose $\limsup_{n \to \infty}\int V(x) |u_{\lambda_n}|^2 \ge \eps >0$. Passing to a subsequence we
may assume $(u_{\lambda_n})_{n \in \nz}$ converges weakly in $E$ to some $u_0 \neq 0$ satisfying
\begin{align*}
\|\nabla u_0\|_2^2 - \lambda_1(\rz^N,V) \int V(x) |u_0|^2 &\le \liminf_{n \to \infty}
I_{\lambda_n}(u_{\lambda_n})\le 0,  
\end{align*}
contradicting $\lambda_1(\rz^N,V)<\lambda_1(\Omega^{-0},V)$ and $\int h(x) |u_0|^p =0$.
Consequently from (\ref{eq:175}) we deduce that $u_{\lambda_n}\to 0$ in $D^{1,2}(\rz^N)\cap
L^p(\rz^N,h)\cap L^2(\rz^N,V^+)$. We use again the fact that $I_\lambda(u_\lambda)<0$ to see that
\[\int V^-(x) |u_{\lambda_n}|^2 \le \int V^+(x)|u_{\lambda_n}|^2 -\frac{1}{\lambda_n}\left(\|\nabla
  u_{\lambda_n}\|_2^2 + \frac{2}{p} \int h(x) |u_{\lambda_n}|^p\right),\]
which finally gives $u_{\lambda_n}\to 0$ in $E$.\\
Fix $(\lambda_n)_{n \in \nz}$ such that $\lambda_n \to \lambda_1(\Omega^{-0},V)$ as $n \to
\infty$. If $\lambda_1(\Omega^{-0},V)=\infty$, then by (\ref{eq:176}) we see
$I_{\lambda_n}(u_{\lambda_n})\to -\infty$ and 
\[-\infty = \lim_{n \to \infty} I_{\lambda_n}(u_{\lambda_n}) - \frac{1}{2}
I'_{\lambda_n}(u_{\lambda_n})u_{\lambda_n} =
\big(\frac{1}{p}-\frac{1}{2}\big) \int h(x) |u_{\lambda_n}|^p.\]  
Consequently $\int h(x) |u_{\lambda_n}|^p \to \infty$ as $n \to
\infty$.\\
Suppose $\lambda_1(\Omega^{-0},V)<\infty$ and $(u_{\lambda_n})_{n \in \nz}$ is bounded in
$E$. Thus by passing to a subsequence we 
may assume that $(u_{\lambda_n})_{n \in \nz}$ converges weakly to $u_0$ in $E$. By (\ref{eq:176})
\[I_{\lambda_1(\Omega^{-0},V)}(u_0) \le \liminf_{n \to \infty} I_{\lambda_n}(u_{\lambda_n})<0.\]    
Thus $u_0\neq 0$ is a nonnegative weak solution of \ref{eq:1} for
$\lambda \ge \lambda_1(\Omega^{-0},V)$, which is impossible due to Theorem \ref{sec:intro:t3} and Lemma
\ref{sec:intro:l1}.
\end{proof}

\begin{lemma}
\label{sec-existence-fred}
Under the assumptions of Theorem \ref{sec:intro:t2} the operator
$I_\lambda''(u_\lambda):E \to E'$ is injective. If, moreover, 
\begin{align}
\label{eq:57}
L^p(\rz^N,h^+) \not \embed D^{1,2}(\rz^N)\cap L^2(\rz^N,|V|), \text{ i.e. } 
D^{1,2}(\rz^N)\cap L^2(\rz^N,|V|) \neq E   
\end{align}
then the inverse of $I_\lambda''(u_\lambda)$, defined on $I_\lambda''(u_\lambda)(E)$, is not continuous.
\end{lemma}
\begin{proof}
To show the injectivity it is enough to prove that 
\begin{align}
\label{eq:60}
(I''_\lambda(u_\lambda)\phi,\phi)>0 \text{ for every } \phi \in E\backslash\{0\}.  
\end{align}
By (\ref{snone1}) applied to $k(x):=\lambda V(x)-h(x)u_\lambda^{p-2}$ we obtain
\begin{align*}
F(\phi) &:= \int |\nabla \phi|^2 - \lambda \int V(x) \phi^2 + \int h(x)u_\lambda^{p-2}(x)\phi^2
\ge 0 \qquad \forall \phi \in E.   
\end{align*}
We may write
\begin{align*}
(I''_\lambda(u_\lambda)\phi,\phi) &= F(\phi)+ (p-2)\int h(x)u_\lambda^{p-2}(x)\phi^2.    
\end{align*}
Thus $(I''_\lambda(u_\lambda)\phi,\phi)>0$ for every $\phi \not\in D^{1,2}(\Omega^{-0})$. 
For $\phi \in D^{1,2}(\Omega^{-0})$ we obtain
\begin{align*}
(I''_\lambda(u_\lambda)\phi,\phi) &= \int |\nabla \phi|^2 - \lambda \int V(x) \phi^2 \ge 
(\lambda_1(\Omega^{-0},V)-\lambda) \int |\nabla \phi|^2,
\end{align*}
and (\ref{eq:60}) follows.\\
To prove the second part of our claim we take a sequence $(\phi_n)_{n \in \nz}\subset E$ such
that $\|\phi_n\|_E=1$ for all $n \in \nz$ and 
$(u_n)_{n \in \nz}$ converges to zero strongly in $D^{1,2}(\rz^N)\cap L^2(\rz^N,|V|)$. This is
possible by (\ref{eq:57}). 
We obtain using H\"older's inequality
\begin{align*}
\sup_{\|\psi\|_E=1} \int h(x) u_\lambda^{p-2}\phi_n \psi 
&\le \sup_{\|\psi\|_E=1} \big( \int h(x)
u_\lambda^{\frac{(p-2)p}{p-1}}|\phi_n|^{\frac{p}{p-1}}\big)^{\frac{p-1}{p}}
\big( \int h(x)|\psi|^p\big)^{\frac{1}{p}}\\
&\le \big( \int h(x)
u_\lambda^{\frac{(p-2)p}{p-1}}|\phi_n|^{\frac{p}{p-1}}\big)^{\frac{p-1}{p}}.   
\end{align*}
Obviously 
$(\phi_n)_{n \in \nz}$ converges to zero pointwise almost everywhere and we may use, similar to
the proof of Lemma \ref{sec:growth-l2}, Vitali's
convergence theorem \cite{HewStr75} to deduce that the latter integral tends to
zero as $n \to \infty$. Consequently, $I_\lambda''(u_\lambda)\phi_n$ tends to zero in $E'$,
which yields the claim.
\end{proof}

\begin{remark}
\label{sec-existence-rem1}
Under assumption (\ref{eq:57}) the operator $I_\lambda''(u_\lambda)$ is not invertible on its
range, hence we cannot use the implicit function theorem to prove that the map $\lambda \mapsto
u_\lambda$ is $C^1(J,E)$.\\ 
Moreover, we may deduce that $I_\lambda''(u_\lambda)$ cannot be a Fredholm operator under
assumption (\ref{eq:57}). Conversely, suppose that the codimension of the range of
$I_\lambda''(u_\lambda)$ is finite. Then we could decompose $E'$ topologically into
$E'=I_\lambda''(u_\lambda)(E)\oplus \mathcal{F}$. 
Since $I_\lambda''(u_\lambda)$ is injective, the open mapping theorem would yield in this case
a continuous inverse on $I_\lambda''(u_\lambda)(E)$, in contradiction to Lemma
\ref{sec-existence-fred}.       
\end{remark}

\begin{proof}[\bf{Proof of Theorem \ref{sec:intro:t2:wachs}}]
Because $V$ is nonnegative we see that $u_{\lambda'}$ is a strict super-solution of \ref{eq:1} for any
$\lambda'> \lambda$. Hence by Lemma \ref{sec:growth-l2} there holds $u_{\lambda}> u_{\mu}$
for any $\lambda> \mu$, which gives 
\begin{align}
\label{eq:187}
v_\lambda(x):=\frac{u_\lambda(x)-u_\mu(x) }{\lambda-\mu} > 0 \text{ for all } x \in \rz^N.  
\end{align}
To estimate the growth behavior of $u_\lambda$ we use an idea given in
\cite{Logo00,GarLogoetal98} and try to differentiate \ref{eq:1} with 
respect to $\lambda$. Because we do not know that the map $\lambda \mapsto u_\lambda$ is $C^1(J,E)$,
we have to work with $v_\lambda$. We obtain for all $\varphi \in E$ and $\lambda>\mu$  
\begin{align}
\label{eq:186}
\begin{split}
\left(T(v_\lambda),\varphi\right) &:=
\int \nabla v_\lambda \nabla\varphi + \int \left(- \lambda V(x) +
 h(x) f_{\lambda,\mu}(x)\right)v_\lambda\varphi= \int V(x) u_\mu \varphi,  
\end{split}
\end{align}
where $f_{\lambda,\mu}(x) := \left(u_\lambda^{p-1}(x)-u_\mu^{p-1}(x)\right)(u_\lambda(x)-u_\mu(x))^{-1}$.\\
Fix $\mu \in (\lambda_1(\rz^N,V),\lambda_1(\Omega^{-0},V))$ and a principal
eigenfunction $e_1(\Omega^{-0},V)$. Because $u_\mu$
is a positive function and $\Omega^{-0}$ is bounded, there is a $C>0$
such that $u_\mu \ge C e_1(\Omega^{-0},V)$.  
Consider the function $w_\lambda\in E$ defined by
\begin{align}
\label{eq:189}
w_\lambda:= v_\lambda- \frac{C}{\lambda_1(\Omega^{-0},V)-\lambda}
e_1(\Omega^{-0},V).  
\end{align}
Consequently, by (\ref{eq:186}), we have for all nonnegative $\varphi \in
D^{1,2}(\Omega^{-0})\cap L^2(\rz^N,V)$  
\begin{align}
\label{eq:188}
\left(T(w_\lambda),\varphi\right)
= \int \left(V(x)u_\mu - C V(x) e_1(\Omega^{-0},V)\right)\varphi \ge 0.   
\end{align}
Testing (\ref{eq:188}) with $w_\lambda^-= -\min(w_\lambda,0) \in
D^{1,2}(\Omega^{-0})\cap L^2(\rz^N,V)$ we see 
\[0\le \left(T(w_\lambda),w_\lambda^-\right)=-\left(T(w_\lambda^-),w_\lambda^-\right) \le 0,\]
because $\lambda<\lambda_1(\Omega^{-0},V)$. Hence $w_\lambda \ge 0$ and the claim follows.  
\end{proof}


\section{Results for sign-changing $h$}
\label{sec:Results-sign-chang}
We give existence and multiplicity results concerning sign-changing $h$. To this end we need to
introduce the set $\Omega^0$, defined by
\[\Omega^0 := \big\{ x \in \rz^N \where h(x) =0\big\}.\]
\begin{lemma}
\label{s>0p1}
Suppose (\ref{eq:2})-(\ref{eq:4}), $h^- \not\equiv 0$ and $\lambda_1(\rz^N,V)<\lambda_1(\rz^N,V,h)<\infty$.
Then we may choose $\alpha >0$ such that $\alpha \cdot u_0$ is a positive solution of
\ref{eq:1} for $\lambda= \lambda_1(\rz^N,V,h)$, where $u_0$ is the minimizer of
$\lambda_1(\rz^N,V,h)$ obtained in Lemma \ref{s>0l1}.
\end{lemma}
\begin{proof}
First we shall show that $\lambda_1(\rz^N,V,h)<\lambda_1(\Omega^0,V)$. Conversely, suppose that this
is wrong. (\ref{eq:4}) implies that $\lambda_1(\Omega^{0},V)$ is attained by some $e_1(\Omega^0,V)$. Because 
$\lambda_1(\rz^N,V)<\lambda_1(\rz^N,V,h)$ we may choose $\varphi_1 \in E$ such that 
\begin{align}
\label{eq:192}
\int \nabla e_1(\Omega^0,V) \nabla \varphi_1 - \lambda_1(\Omega^0,V) \int V(x) e_1(\Omega^0,V)
\varphi_1 <0.  
\end{align}
Because $h^- \not\equiv 0$ there is $\varphi_2 \in C_c(\{x \where h(x)<0\})$ such that
\begin{align}
\label{eq:193}
\int h(x) |\varphi_1+\varphi_2|^p <0.  
\end{align}
By (\ref{eq:192}) and (\ref{eq:193}) we see for small values of $t>0$
\begin{align*}
\int h(x)\big|e_1(\Omega^{0},V)+t(\varphi_1+\varphi_2)\big|^p=t^p \int
h(x) |\varphi_1+\varphi_2|^p <0  \text{  and }\\ 
\frac{\|\nabla(e_1(\Omega^{0},V)+t(\varphi_1+\varphi_2))\|_2^2}{\int
  V(x)(e_1(\Omega^{0},V)+t(\varphi_1+\varphi_2))^2}
< \frac{\|\nabla e_1(\Omega^{0},V)\|_2^2}{\int V(x) e_1(\Omega^{0},V)} = \lambda_1(\Omega^0,V).
\end{align*}
Thus $\lambda_1(\rz^N,V,h)< \lambda_1(\Omega^0,V)$.\\
To obtain a contradiction, suppose ${\int h(x) |u_0|^{p}<0}$. Because the
set $\{u \in E\where \int h(x) |u|^{p} <0\}$ is
open in $E$, the function $u_0$ is a local minimizer
of ${\inf\{\|\nabla u\|^{2}_2\where\int V(x) |u|^2 = 1\}}$.
Therefore $u_0$ is a weak solution of
${-\laplace u = \lambda_1(\rz^{N},V,h) V(x) u}$ in $\rz^{N}$. Standard regularity results and Harnack's
inequality (see \cite[C.1]{Sim82}) show that $u_0$ is a positive and continuous function.
Now the estimate (\ref{snone1}) in Lemma \ref{snonl1} contradicts the assumption
${\lambda_1(\rz^N,V,h)>\lambda_1(\rz^N,V)}$.\\
Thus we have ${\int h(x) |u_0|^{p}=0}$. We define $M: E \to \rz^{2}$ by
\[M(u) := (M_1(u), M_2(u)) = \left(\frac{1}{2} \int V(x) |u|^2, \frac{1}{p} \int
h(x) |u|^{p}\right).\]
Obviously $M \in C^{1}(E,\rz^{2})$ and for all $u,\varphi \in E$ there holds
\[M'(u) \varphi = \left(\int V(x) u \varphi, \int h(x) |u|^{p-2}u \varphi\right).\]
We have ${\displaystyle M'(u_0)u_0 = (1,0)}$. If ${\int h(x)
|u_0|^{p-2}u_0 \varphi =0 }$ for all $\varphi \in E$,
then ${\displaystyle h(x) |u|^{p-2}u}$ is zero almost everywhere, which implies
$u_0(x)=0$ for almost every $x \in \rz^{N} \backslash \Omega^{0}$, contradicting
$\lambda_1(\rz^{N},V,h) <\lambda_1(\Omega^{0},V).$
The above arguments show that $M'(u_0)$ is surjective. Because $u_0$ minimizes 
${\inf\left\{\frac{1}{2} \|\nabla u\|_2^{2}\where M(u)=(\frac{1}{2},0)\right\} }$ we deduce from the
Lagrange rule the existence of  
${\displaystyle \sigma= (\sigma_1,\sigma_2) \in \rz^{2}}$ such that for all $\varphi \in E$ 
\begin{align}
\label{eq:44} 
\int \nabla u_0 \nabla \varphi = \sigma_1 \int V(x) u_0 \varphi\, + \sigma_2 \int 
h(x) |u_0|^{p-2} u_0 \varphi. 
\end{align} 
Testing (\ref{eq:44}) with $u_0$ we get $\sigma_1= \lambda_1(\rz^N,V,h)$. Define ${
I_{RQ}(u) := \frac{\|\nabla u\|_2^{2}}{\int V(x) |u|^2}}$. Then $u_0$ minimizes
${\{I_{RQ}(u)\where\int V(x) |u|^2 > 0, \; \int h(x) |u|^{p}\le 0\}}$. Take a
$\varphi \in E$ with ${\int h(x)
|u_0|^{p-2}u_0 \phi <0 }$. Then we get for all small positive $t$ that
$\int h(x) |u_0+t \varphi|^{p} \le 0$.
The function $u_0$ is a minimizer of $I_{RQ}$ in direction of $\varphi$, this implies
\[0 \le \frac{1}{2} I_{RQ}'(u_0)\varphi = \int \nabla u_0 \nabla \varphi - \lambda_1(\rz^{N},V,h)
\int V(x) u_0 \varphi = \sigma_2 \int h(x) |u_0|^{p-2} u_0 \varphi. \]
Hence we see $\sigma_2 \le 0$, which implies $\sigma_2 <0$, because otherwise we get
\[-\laplace u_0 = \lambda_1(\rz^{N},V,h)\, V(x) u_0 \text{ in } \rz^{N},\]
which is impossible as we have seen before.
The claim follows for $\alpha := (-\sigma_2)^{\frac{1}{p-2}}$.
\end{proof}
\begin{corollary}
\label{sec:homogeneous-case:cl}
Under the assumptions of Lemma \ref{s>0p1} the equation \ref{eq:1} admits an ordered pair of
solutions $u_0<u_1$ for $\lambda= \lambda_1(\rz^N,V,h)$ such that 
\[I_{\lambda_1(\rz^N,V,h)}(u_0)<0=I_{\lambda_1(\rz^N,V,h)}(u_1)= 
\big(I''_{\lambda_1(\rz^N,V,h)}(u_1)u_1,u_1\big).\]
\end{corollary}
\begin{proof}
We let $u_1$ be the solution given by Lemma \ref{s>0p1}. The second solution $u_0$ exists due
to Lemma \ref{sec:growth-l2}.
\end{proof}

\section*{Second solution}
We show that $I_{\lambda}$ admits a mountain-pass geometry. To this
end we fix $\lambda>0$ and suppose throughout 
the rest of the section that the assumptions of Theorem
~\ref{sec:intro:t1} are satisfied, i.e. we suppose
(\ref{eq:2})-(\ref{eq:4}) and 
\begin{align}
\label{eq:69}
\lambda_1(\rz^N,V) < \lambda < \lambda_1(\rz^N,V,h) 
\end{align}
Moreover, we assume $h^- \not\equiv 0$, which is a necessary condition for the
existence of a second positive solution to \ref{eq:1} (see Theorem \ref{sec:intro:t1}).
Because $V \in C(\rz^N)$ and $h^-(x_0)\neq 0$ for some $x_0 \in \rz^N$ it is possible to choose
a $C_c^\infty(\rz^N)$-function $\varphi_0$ concentrated at $x_0$ such that
\begin{align}
\label{eq:71}
\lambda_1(\rz^N,V,h) \int V(x) |\varphi_0|^2 < \int |\nabla \varphi_0|^2 \text{ and }
I_{\lambda}(\varphi_0)\le 0.  
\end{align}
We define a set of paths
\[\Gamma := \left\{\gamma \in C^1([0,1],E)\where \gamma(0)= 0,\, 
\gamma(1)=\varphi_0 \right\}.\] 
\begin{lemma}
\label{sec:Second-solution:l1}
\[c:= \inf_{\gamma \in \Gamma} \sup_{t \in [0,1]} I_{\lambda}(\gamma(t)) >0\]  
\end{lemma}
\begin{proof}
Set $M^+:= \left\{u \in E\where \lambda_1(\rz^N,V,h) \int V(x) |u|^2
  \le \|\nabla u\|_2^{2} \right\}\neq \emptyset$.  
First we observe 
\begin{equation}
\label{sec:Second-solution:e1}
\|\nabla u\|_2^2 - \lambda \int V(x) |u|^2 \ge
\left(1-\frac{\lambda}{\lambda_1(\rz^N,V,h)}\right)\|\nabla u\|_2^2  
\qquad \forall \, u \in M^+.
\end{equation}
With the help of (\ref{sec:Second-solution:e1}) we derive for all $u \in M^+$
\begin{align*}
I_{\lambda}(u)&\ge C_1(\lambda)\|\nabla u\|^2_2 - \frac{1}{p} \int h^-(x) |u|^p 
\ge  C_1(\lambda) \|\nabla u\|^2_2 +
o(\|\nabla u\|^2_2)_{\|\nabla u\|_2\to 0}.   
\end{align*}
Thus there exists $r>0$ such that $\|\nabla\varphi_0\|_2>r>0$ and
\begin{align}
\label{eq:70}
 I_{\lambda}(u)\ge c_1 >0 \mbox{  for all } u \in M^+  \text{ such
 that }\|\nabla u\|_2 =r.
\end{align}
Suppose there is a $u_1 \in E$ such that ${\displaystyle \|\nabla u_1\|_2^2 = \lambda_1(\rz^N,V,h) \int V(x) |u_1|^2}$. 
Then the definition of $\lambda_1(\rz^N,V,h)$ and (\ref{sec:Second-solution:e1}) imply
\begin{align}
\label{sec:Second-solution:e2}
\begin{split}
I_{\lambda}(u_1) &= \frac{1}{2}\left(\|\nabla u_1\|_2^2-\lambda\int V(x) |u_1|^2\right)
+\underbrace{\frac{1}{p}\int h(x) |u_1|^p}_{\ge 0}\\
&\ge \frac{1}{2}\left(\|\nabla u_1\|_2^2-\lambda\int V(x) |u_1|^2\right)
\ge \frac{1}{2}\left(1-\frac{\lambda}{\lambda_1(\rz^N,V,h)}\right) \|\nabla u\|^2_2   
\end{split}
\end{align}
Let $\gamma$ be an arbitrary element of $\Gamma$. Let
\[t_0 := \max\left\{t \in [0,1] \where \|\nabla(\gamma(t))\|_2^2 = 
\lambda_1(\rz^N,V,h) \int V(x) |\gamma(t)|^2\right\}.\]
Our choice of $\varphi_0$ yields that $t_0$ is attained and $0\le t_0<1$.
If $\|\gamma(t_0)\|_{D^{1,2}(\rz^N)}\ge r$, then (\ref{sec:Second-solution:e2}) shows 
\[\sup_{t \in [0,1]}I_{\lambda}(\gamma(t)) \ge
\frac{1}{2}\left(1-\frac{\lambda}{\lambda_1(\rz^N,V,h)}\right) r^2
>0.\] 
If $\|\gamma(t_0)\|_{D^{1,2}(\rz^N)}< r$, then $\gamma(t) \in M^+$ for all $t \in [t_0,1]$ 
and (\ref{eq:70}) shows 
\[\sup_{t \in [0,1]}I_{\lambda}(\gamma(t)) \ge c_1 >0.\]
Thus the claim follows. 
\end{proof}

By Lemma \ref{sec:Second-solution:l1} the functional $I_{\lambda}$ has
a mountain pass geometry. Consequently using versions of the classical
mountain pass theorem of \citet{AmbRab73} given in \cite{Cer78,BarBenFor83} there exits a $(PS)_c$ 
sequence $(u_n)_{n \in \nz}$, i.e.
\[I_{\lambda}(u_n) \xpfeil{n \to \infty} c \text{ and }
I_{\lambda}'(u_n)(1+\|u_n\|_E) \xpfeil{n \to \infty} 0.\] 
Before we show that such $(PS)_c$ sequences give rise to critical points of
$I_\lambda$, we observe that it is
sufficient to consider special $(PS)_c$ sequences, consisting of almost nonnegative functions.
\begin{lemma}
\label{sec:second:nonneg}
There exists a $(PS)_c$ sequence $(u_n)_{n \in \nz}\subset E$, i.e.
\[ I_{\lambda}(u_n) \xpfeil{n \to \infty} c>0,\; I_{\lambda}'(u_n)(1+\|u_n\|_E)
\xpfeil{n \to \infty} 0 \text{ in }E',\] 
such that $\dist(u_n, E^+) \to 0$ as $n \to \infty$, where we denote by $E^+$ the cone of
nonnegative functions in $E$. In the sequel we will call a sequence with
these properties a $(PS)_c^+$ sequence. 
\end{lemma}
\begin{proof}
Lemma \ref{sec:Second-solution:l1} shows 
\[c:= \inf_{\gamma \in \Gamma} \sup_{t \in [0,1]} I_{\lambda}(\gamma(t)) >0.\]
Hence there is a sequence of paths $(\gamma_n)_{n \in \nz}$ such that
\[\lim_{n \to \infty} \max_{u \in \gamma_n} I_{\lambda}(u) =c.\]
Because $I_{\lambda}(u) = I_{\lambda}(|u|)$ we may assume $\gamma(\cdot)=|\gamma(\cdot)|$. A
version of the classical mountain pass theorem in 
\cite[Th. 13.10]{KuzPoh97} shows that there is a $(PS)_c$ sequence $(u_n)_{n \in \nz}$ such that 
$\|u_n - \gamma_n\|_{E} \to 0$ by using a deformation argument in the vicinity of almost
minimizing paths, which carries over to the work of \cite{Cer78,BarBenFor83}. Hence the claim follows.  
\end{proof}

\begin{lemma}
\label{sec:Second-solution-2l}
Suppose $(u_n)_{n \in \nz}\subset E$ is a $(PS)^+_c$ sequence.
Then there exists a nonnegative $u_0 \in E$, such that
$I_{\lambda}(u_0)=c$, $I_{\lambda}'(u_0)=0$ and we may assume  
after going to a subsequence that $(u_n)_{n \in \nz}$ converges in $E$ to $u_0$.  
\end{lemma}
\begin{proof}
Let $(u_n)_{n \in \nz}$ be a $(PS)^+_c$ sequence as stated above.  
We have
\begin{align}
\label{eq:34}
-2c &= \lim_{n \to \infty}
\left(I'_\lambda(u_n)u_n-2I_\lambda(u_n)\right) 
= \lim_{n \to \infty} \frac{p-2}{p} \int h(x) |u_n|^p.  
\end{align}
Hence the definition of $\lambda_1(\rz^N,V,h)$ yields
\begin{align}
\label{eq:45}
c &= \lim_{n \to \infty}
\left(I_\lambda(u_n)-\frac{1}{p}I_\lambda(u_n)u_n\right) \ge
\frac{p-2}{2p}\left(1-\frac{\lambda}{\lambda_1(\rz^N,V,h)}\right) 
\limsup_{n\to\infty} \|\nabla u_n\|_2^2   
\end{align}
Consequently $(u_n)_{n \in \nz}$ is bounded in $D^{1,2}(\rz^N)$ and by
(\ref{eq:3}) in $L^p(\rz^N,h^-)$. The estimate (\ref{eq:34}) now implies that
$(u_n)_{n \in \nz}$ is bounded in $L^p(\rz^N,h^+)$. Finally
(\ref{eq:4}) and the boundedness of $(I_\lambda(u_n))_{n \in \nz}$ show
that $(u_n)_{n \in \nz}$ is bounded in $E$.\\  
Because $(u_n)_{n \in \nz}$ is bounded in $E$, we may choose a subsequence (still denoted by
$(u_n)_{n \in \nz}$), which converges weakly to $u_0$ in $E$ and pointwise almost
everywhere. Since $(u_n)_{n \in \nz}$ is a $(PS)^+_c$ sequence we have
$\|u_n-v_n\|_{E}\to 0$ for some nonnegative $v_n \in E$   
and $v_n \schwach u_0$ as $n \to \infty$. Thus $u_0$ is nonnegative, for the cone of
nonnegative functions is weakly closed.
The sequence $(u_n)_{n \in \nz}$ converges weakly in $E$ and pointwise
almost everywhere to $u_0 \in E$. Consequently, using for instance
Vitali's convergence theorem \cite{HewStr75}, we have $I'_\lambda(u_0)=0$.
Moreover
\begin{align}
o(1) &= \left(I_{\lambda}'(u_n)-I_{\lambda}'(u_0)\right)(u_n-u_0)\nonumber \\
&= \|\nabla (u_n-u_0)\|_2^2 + \lambda \int V^-(x) |u_n-u_0|^2  
- \lambda \underbrace{\int V^+(x) |u_n-u_0|^2}_{\to 0\; (n\to \infty) \text{ by (\ref{eq:4})}}\nonumber\\
&\quad -\underbrace{\int h^-(x)(|u_n|^{p-2}u_n-u_0^{p-1})(u_n-u_0)}_{\to
  0\; (n\to \infty) \text{ by (\ref{eq:3})} } 
 + \int h^+(x)(|u_n|^{p-2}u_n-u_0^{p-1})(u_n-u_0)\nonumber\\
&= o(1)+ \|\nabla (u_n\!-\!u_0)\|_2^2+ \lambda \int V^-(x) |u_n\!-\!u_0|^2\!+\!\int\!h^+(x)
\underbrace{(|u_n|^{p-2}u_n\!-\!u_0^{p-1})(u_n\!-\!u_0)}_{\ge 0} \nonumber\\
\label{eq:75}
\end{align}
This implies ${\displaystyle u_n \to u_0}$ in $D^{1,2}(\rz^N)\cap
L^2(\rz^N,|V|)$.
By (\ref{eq:75}) and the Brezis-Lieb lemma \cite{BreLie83} we obtain
\[0 =\lim_{n \to \infty} \int h^+(x)(|u_n|^{p-2}u_n\!-\!u_0^{p-1})(u_n\!-\!u_0) = 
 \lim_{n \to \infty}\left(\int h^+(x) u_n^{p} - \int h^+(x) u_0^{p}\right),\]
and the uniform convexity of $L^p$ implies $u_n \to u_0$ in $L^p(\rz^N,h^+)$.
Hence we finally see $(u_n)_{n \in \nz}$ converges to $u_0$ in $E$.   
\end{proof}

\begin{proof}[\bf{Proof of Theorem \ref{sec:intro:t3}}]
The assertion of Theorem \ref{sec:intro:t3} for
$\lambda=\lambda_1(\rz^N,V,h)$ is an immediate consequence of Corollary
\ref{sec:homogeneous-case:cl}. Note that any solution given by
Lemma \ref{sec:growth-l2} satisfies
\[(I_\lambda''(u_\lambda)\phi,\phi)\ge 0 \text{ for all } \phi \in
E.\]
For $\lambda_1(\rz^N,V)<\lambda<\lambda_1(\rz^N,V,h)$ Lemma
\ref{sec:Second-solution-2l} yields the existence of $v_\lambda$.
Moreover, any solution $v_\lambda$, such that $I_\lambda(v_\lambda)>0$, satisfies due to the homogeneity of
the nonlinearity
\[(I_\lambda''(v_\lambda)v_\lambda,v_\lambda)<0.\]
The claim now follows using again Lemma \ref{sec:growth-l2}.   
\end{proof}

\section{Existence for $\lambda=\lambda^*$}
\label{sec:lambda}
We show that \ref{eq:1} is solvable for $\lambda=\lambda^*$ by approximating the desired solution
with solutions for $\lambda<\lambda^*$. To justify the passage to the limit we will need that
$\lambda^*<\lambda_1(\Omega^0,V)$, where $\Omega^0=\{x \in \rz^N \where h(x)=0\}$.\\    
If $h$ is nonnegative, then $\lambda^*$ equals $\lambda_1(\Omega^0,V)$. If $h^- \neq 0$ then we
have $\Omega^0 \subsetneqq \Omega^{-0}$, hence $\lambda_1(\Omega^0,V)>\lambda_1(\Omega^{-0},V)\ge \lambda^*$ by
the strict monotonicity of the first eigenvalue with respect to the domain. This is a nice, but
unfortunately wrong argument because $\Omega^{-0}$ need not to be a domain. A counterexample may 
easily be obtained by choosing a function $h$ vanishing on a large ball $B_1$ and changing sign only on a
smaller second ball $B_2$ disjoint from $B_1$. In that case the eigenfunction corresponding to the
first eigenvalue vanishes on the set where $h$ attains negative values and $\lambda_1(\Omega^{-0},V)$
equals $\lambda_1(\Omega^0,V)$. The next lemma gives some sufficient conditions ensuring that
$\lambda^*<\lambda_1(\Omega^0,V)$, although we believe that 
$\lambda^*<\lambda_1(\Omega^0,V)$ whenever $h^-\not\equiv 0$.  
\begin{lemma}
\label{sec:results-sign-chang:l2}
Under the assumptions of Theorem \ref{sec:intro:t1}, if moreover $h^- \not\equiv 0$ and either
\begin{align}
\label{eq:190}
\rz^N\backslash \supp(h^+) \text{ is connected or}\\
\label{eq:191}
V\ge 0 \text{ and } 0<\int h^-(x) e_1(\Omega^{-0},V)^p
\end{align}
for some principal eigenfunction $e_1(\Omega^{-0},V)$. Then $\lambda^* < \lambda_1(\Omega^0,V)$.
\end{lemma}
\begin{proof}
Under assumption (\ref{eq:190}) Harnack's inequality shows that any principal eigenfunction
$e_1(\Omega^{-0},V)$ is positive on $\{x\in \rz^N \where h(x)<0\}$. Consequently
$\lambda_1(\Omega^{-0},V)<\lambda(\Omega^0,V)$. From Theorem \ref{sec:intro:t3} we infer
$\lambda^*\le \lambda_1(\Omega^{-0},V)$.\\ 
Under the hypothesis (\ref{eq:191}) we suppose contrary to the
assertion of the lemma that there is a sequence  
$(u_n,\lambda_n)_{n \in \nz}$ in $E \times (\lambda_1(\rz^N,V),\lambda_1(\Omega^0,V))$ such that
$\lambda_n \to \lambda_1(\Omega^0,V)$ monotone increasing as $n \to
\infty$ and $u_n$ is a solution of \ref{eq:1} for 
$\lambda= \lambda_n$. Lemma \ref{sec:intro:l1} and Theorem \ref{sec:intro:t3} show $\lambda_n <
\lambda_1(\Omega^{-0},V)$. Hence $\lambda_1(\Omega^{-0},V)=\lambda_1(\Omega^{0},V)$. By Lemma
\ref{sec:growth-l2} we may assume  
\[\left(I_{\lambda_n}''(u_n)\varphi,\varphi\right)\ge 0 \quad \forall \varphi \in E \text{ and }
I_{\lambda_n}(u_n)<0.\]  
Because $V$ is nonnegative, each $u_n$ is a super-solution of \ref{eq:1} with $\lambda\le \lambda_n$
and $h$ replaced by $h^+$.
We may use Lemma \ref{sec:growth-l2} and Theorem \ref{sec:intro:t2}
to derive for all $n$ that $u_n \ge u_0$, where $u_0$ is the unique positive solution of \ref{eq:1}
with $\lambda= \lambda_2$ and $h$ replaced by $h^+$. Consequently 
\begin{align*}
0 &\le \liminf_{n \to \infty}
\left(I_{\lambda_n}''(u_n)e_1(\Omega^{-0},V),e_1(\Omega^{-0},V)\right)
= -(p-1) \liminf_{n \to \infty} \int h^-(x) |u_n|^{p-2}e_1(\Omega^{-0},V)^2\\
&\le -(p-1) \int h^-(x) |u_0|^{p-2} e_1(\Omega^{-0},V)^2,  
\end{align*}
contradicting (\ref{eq:191}).
\end{proof}

\begin{lemma}
\label{sec:results-sign-chang:l1}
Suppose (\ref{eq:2})-(\ref{eq:4}), $h^- \not\equiv 0$ and
$\lambda_1(\rz^N,V)< \lambda^* <\lambda_1(\Omega^0,V)$. Then \ref{eq:1} is solvable for $\lambda=\lambda^*$. 
\end{lemma}
\begin{proof}
Let $(u_n,\lambda_n)_{n \in \nz}$ be a sequence in $E \times ]\lambda_1(\rz^N,V),\lambda^*]$ such
that $\lambda_n \to \lambda^*$ as $n \to \infty$ and $u_n$ is a positive solution of \ref{eq:1} for
$\lambda= \lambda_n$. By Lemma \ref{sec:growth-l2} we may assume 
\begin{align}
\label{eq:180}
I_{\lambda_n}(u_n)<0 \text{ and } (I_{\lambda_n}''(u_n),u_n)\ge 0 \text{ for all } n \in \nz. 
\end{align}
Moreover, we have
\begin{align}
\label{eq:194}
0 \le (I_{\lambda_n}''(u_n),u_n)-2I(u_n) \le (p-2) \int h(x) |u_n|^p.   
\end{align}
From (\ref{eq:180}) and (\ref{eq:194}) we deduce that $(u_n)_{n \in \nz}$ is bounded in $E$ if
$\left(\int V(x) u_n^2 \right)$ is bounded above. To obtain a contradiction we suppose $\sigma_n^2:=
\int V(x)u_n^2 \to \infty$ as $n \to \infty$ and write $w_n:= u_n/\sigma_n$. Obviously $(w_n)_{n \in
  \nz}$ is bounded in $E$ and we may assume $(w_n)_{n \in \nz}$ converges weakly to some $w_0\neq 0$ in
$E$. There holds for all $\varphi \in E$
\begin{align}
\label{eq:195}
0= \frac{1}{\sigma_n^2}I'(u_n)\varphi = \int \nabla w_n \nabla \varphi - \lambda \int V(x) w_n
\varphi + \sigma_n^{p-2} \int h(x) w_n^{p-1}\varphi.  
\end{align}
Hence $w_0(x) =0$ for all $x$ such that $h(x) \neq 0$, because we have for all $\varphi \in E$
\[\int h(x) w_0(x)^{p-1} \varphi = \lim_{n \to \infty} \int h(x) w_n^{p-1}\varphi =0.\]   
By (\ref{eq:180}) and (\ref{eq:194}) we see $\|\nabla w_0\|_2^2 - \lambda^* \int V(x) w_0^2 \le 0$,
contradicting $\lambda^*<\lambda_1(\Omega^0,V)$.\\
Consequently $(u_n)_{n \in \nz}$ is a bounded $(PS)^+_c$ sequence for $I_{\lambda^*}$ and some $c
\le 0$. By Lemma \ref{sec:Second-solution-2l} and passing to a subsequence we may assume
$(u_n)_{n\in \nz}$ converges to some $u_0$ in $E$. Because each $u_n$ is a super-solution of
\ref{eq:1} with $\lambda\le\lambda_n$ and $h$ replaced by $h^+$, we see that $u_n \ge w_0$, where
$w_{0}$ denotes the unique positive solution of \ref{eq:1} with $\lambda=\lambda_2$ and $h$
replaced by $h^+$. Consequently $u_0=0$ is 
impossible and the claim follows.  
\end{proof}

\begin{proof}[\bf{Proof of Theorem \ref{sec:intro:t4}}]
Lemma \ref{sec:results-sign-chang:l2} shows that the assumptions of Lemma
\ref{sec:results-sign-chang:l1} are satisfied. Hence \ref{eq:1} admits a solution $u_{\lambda^*}$ for
$\lambda=\lambda^*$. Because $V$ is nonnegative $u_{\lambda^*}$ is a super-solution to \ref{eq:1} for
every $\lambda \le \lambda^*$. Consequently Lemma \ref{sec:growth-l2} yields the validity of
(i). The assertion of (ii) immediately follows from Theorem \ref{sec:intro:t5} and the definition of
$\lambda^*$.
\end{proof}

\section*{Acknowledgements}
I would like to thank Prof. Hans-Peter Heinz for many helpful suggestions during the
preparation of the paper.

\bibliographystyle{adinat}
\bibliography{calc_indef}

\end{document}